\documentclass{article}
\usepackage[a4paper, margin=0.5in]{geometry}
\usepackage{mathpazo,mathrsfs, graphicx, framed, color, amssymb, amsmath,amsthm, amsbsy,upgreek}
\usepackage[boxed, noline, ruled, linesnumbered]{algorithm2e}
\usepackage[colorlinks, citecolor=blue]{hyperref}
\usepackage{diagbox}[2011/11/22]
\usepackage{tabularx}

\newtheorem{remark}{Remark}[section]

\def\cE{\mathcal{E}}

\def\cH{\mathcal{H}}
\def\cJ{\mathcal{J}}

\def\cQ{\mathcal{Q}}
\def\cT{\mathcal{T}}

\def\cR{\mathcal{R}}
\def\cM{\mathcal{M}}

\allowdisplaybreaks 
\newcommand{\tabincell}[2]{\begin{tabular}{@{}#1@{}}#2\end{tabular}}

\begin{document}

\title{On accelerating a multilevel correction adaptive finite element
  method for Kohn-Sham equation\footnote{This research is supported partly by National Key  R\&D Program of China 2019YFA0709600, 2019YFA0709601,  National Natural Science Foundations of China (Grant Nos. 11801021, 11922120, 11871489 and 11771434), MYRG of University of Macau (MYRG2019-00154-FST) and the National Center for Mathematics and Interdisciplinary Science, CAS. }}
\author{Guanghui Hu\footnote{Department of Mathematics, University of Macau, Macao S.A.R., China, and
Zhuhai UM Science \& Technology Research Institute , Zhuhai, Guangdong, China (garyhu@umac.mo).} , \ \ \ \ \
Hehu Xie\footnote{LSEC, ICMSEC, Academy of Mathematics and Systems Science, Chinese Academy of Sciences, Beijing 100190, China, and  School of Mathematical Sciences, University of Chinese Academy of Sciences, Beijing, 100049, China (hhxie@lsec.cc.ac.cn).} \ \ \ \ \
and \ \ \ \ \
Fei Xu\footnote{Faculty of Science, Beijing University of
    Technology, Beijing 100124, China (xufei@lsec.cc.ac.cn).}}

\date{}

\maketitle

\begin{abstract}
  Based on the numerical method proposed in [G. Hu, X. Xie, F. Xu,
    J. Comput. Phys., 355 (2018), 436--449.] for Kohn-Sham equation,
  further improvement on the efficiency is obtained in this paper by
  i). designing a numerical method with the strategy of separately
  handling the nonlinear Hartree potential and exchange-correlation
  potential, and ii). parallelizing the algorithm in an eigenpairwise
  approach. The feasibility of two approaches are analyzed in detail,
  and the new algorithm is described completely. Compared with
  previous results, a significant improvement of numerical efficiency
  can be observed from plenty of numerical experiments, which make the
  new method more suitable for the practical problems.

\vskip0.3cm {\bf Keywords.} Kohn-Sham equation, multilevel correction method, finite element method, separately handing nonlinear terms, parallel computing.
\vskip0.2cm {\bf AMS subject classifications.} 65N30, 65N25, 65L15, 65B99.
\end{abstract}

\section{Introduction}

Kohn-Sham Density functional theory (DFT) is one of the most
successful approximate models in the study of many-body system. Its
application has covered many practical application areas , ranging
from chemistry, physics, materials science, chemical engineering,
etc. With rapid development of the hardware, new algorithms and
acceleration techniques are desired to keep improving the efficiency
of the simulations in density functional theory.

So far, lots of numerical methods for solving Kohn-Sham equation have
been developed. For instance, plane-wave method is
the most popular method in the computational quantum chemistry
community. Owing to the independence of the basis function to the
ionic position, plane-wave method has advantage on calculating
intermolecular force. Combined with the pseudopotential method,
plane-wave method plays an important role in the study of the ground
and excited states calculations, and geometry optimization of the
electronic structures.  Although the plane-wave method is popular
in the computational quantum chemistry community, it is
inefficient in the treatment of non-periodic systems like molecules,
nano-clusters, etc., or materials systems with defects, where higher
basis resolution is often required in some spatial regions and a
coarser resolution suffices elsewhere. Furthermore, the plane-wave
method uses the global basis which significantly affect the
scalability of computations on parallel computing platforms.  The
atomic-orbital-type basis sets
\cite{HehreStewartPople,Jensen,WillsCooper} are also widely used for
simulating materials systems such as molecules and clusters. However,
they are well suited only for isolated systems with special boundary
conditions. It is difficult to develop a systematic basis-set for all
materials systems.  Thus over the past decade, more and more attention
have been attracted to develop efficient and scalable real-space techniques
for electronic structure calculations. For more information, we refer
to
\cite{Beck,BowlerChoudhuryGillanMiyazaki,CastroAppelOliveiraRozziAndradeLorenzen,ChelikowskyTroullierSaad,
  GenoveseVideauOspiciDeutschGodeckerMehaut,
ModineZumbachKaxiras,SkylarisHaynesMostofiPayne,SolerArtachoGaleGarciaJunqueraOrdejnPortal} and references
therein for a comprehensive overview.

Among all those real-space methods for Kohn-Sham equation, the finite
element method is a very competitive one. The advantages of finite
element method for solving the partial differential equations include
that it can use unstructured meshes and local basis sets, and it is
scalable on parallel computing platforms.  So far, the application of
the finite element method on solving Kohn-Sham equation has been
studied systematically. Please refer to
\cite{BaoHuLiu,BylaskaHostWeare,FangGaoZhou,HermannsonYevick,LehtovaaraHavuPuska,LinLuE,MasudKannan,
  PaskKleinFondSterne,PaskKleinSterneFong,PaskSterne,
  SchauerLinder,SuryanarayanaGaviniBlesgenBhattacharyaOrtiz,TsuchidaTsukada_1995,TsuchidaTsukada_1996,TsuchidaTsukada_1998,WhiteWilkinsTeter,
  ZhangShenZhou}, and references therein.

Efficiency improvement is a key issue in the development of the method towards the practical simulations. Many
acceleration techniques such as adaptive mesh methods, multigrid preconditioning in solving the eigenvalue problems,
parallelization, have been studied in depth for the purpose. It is noted that
Xie and his co-workers proposed and developed a multilevel
correction technique for solving eigenvalue problems
\cite{ChenXieXu,JiaXieXu,LinXie,Xie_JCP,Xie_NonlinerEig,XieXie}.  With
this technique, solving the nonlinear eigenvalue problem defined on
the finest mesh can be transformed to solving the linear boundary
value problem on the finest mesh and a fixed and low dimensional
nonlinear eigenvalue problem. Fruitful results have been obtained to
successfully show the capability of the multilevel correction method
on improving the efficiency. In \cite{HuXieXu}, a numerical framework
consisting of the multilevel correction method and the $h$-adaptive
mesh method is proposed for solving the Kohn-Sham equation. Similar to
original idea of the multilevel correction method, the nonlinear
eigenvalue problem derived from the finite element discretization of
the Kohn-Sham equation is fixed in a relatively coarse mesh, while the
finite element space built on this coarse mesh is kept enriching by
the solutions from a series of boundary value problems derived from
the Kohn-Sham equation. Here the $h$-adaptive mesh method is used to
tailor a fitting nonuniform mesh for the derived boundary value
problem, so that the finite element space for the nonlinear eigenvalue
problem can be improved well. The performance of proposed solver for
the Kohn-Sham equation has been checked by the following works
\cite{ChenXieXu,HuXieXu,JiaXieXu,Xie_NonlinerEig,XieXie}.


In this paper, the efficiency of the proposed numerical method in
\cite{HuXieXu} will be further improved, based on following two
approaches. The first approach is based on an observation of
significant difference of contributions for total energy of the system
from two nonlinear terms in the hamiltonian, i.e., Hartree potential
and exchange-correlation potential. It is the nonlinearity introduced
by these two potentials which makes the analysis and calculation of
Kohn-Sham equation nontrivial. Unlike the Hartreen potential which has
an exact expression, there is no exact expression for the
exchange-correlation potential. Hence, an approximation such as local
density approximation (LDA) is needed for the study. An interesting
observation for the difference between Hartree and
exchange-correlation potentials can be made, based on following
results from NIST standard reference database 141
\cite{refdata}.
\begin{table}[!htbp]
  \begin{center}
    \begin{tabularx}{0.6\textwidth}{|X|X|X|X|X|}
      \hline
      E\_tot & E\_kin & E\_har & E\_coul & E\_xc\\
      \hline
      -847.277 & 846.051 & 355.232 & -2008.741 & -39.319\\
      \hline
    \end{tabularx}
    \caption{Ground state total energy and its components of a
      Titanium atom. The data is from
      \cite{refdata}.}
  \end{center}
\end{table}
It is clearly seen that the magnitude of the exchange-correlation
energy (E\_xc) is just around 10\% of the Hartree energy (E\_har). It also should be
noted that similar comparisons exist for all other elements in the
period table. Such observation brings us a chance to further improve
the efficiency of the algorithm by separately considering two
nonlinear terms. The idea is to reduce the computational resource for
handling the nonlinearity introduced by the Hartree
potential. which is realized by designing a new iteration
scheme to replace the traditional self-consistent field iteration scheme in this paper. The new scheme consists of two nested iteration
schemes, i.e., an outer iteration for resolving the nonlinearity from
the exchange-correlation potential, and an inner iteration for
resolving the nonlinearity from the Hartree potential. Although there
are two iterations in our algorithm, it is noted that the outer
iteration process for the nonlinearity of the exchange-correlation
term always be done around 10 times in our numerical experiments for
molecules from a lithium hydride (2 atoms) to a sodium cluster (91
atoms). For the inner iteration, dozens of iterations are needed in
each outer iteration at the beginning. However, with the increment of
the outer iterations, the number of the inner iterations decreases
significantly. Hence, the total number of the iterations (inner
iterations by outer iterations) of our method is fairly comparable
with that of a standard SCF iteration. However, due to the separation
of two nonlinear terms and the framework of multilevel correction
method, the calculation involving the basis function defined on the
fine mesh can be extracted separately and precalculated, so that the
computational work of the inner iteration only depends on the
dimension of the coarse space. With this strategy, a large amount of
CPU time can be saved in the simulations, compared with the original
multilevel correction method \cite{HuXieXu}.

To further improve the efficiency, the eigenpairwise parallelization of the
algorithm based on message passing interface (MPI) is also studied in this work.
One desired feature for the algorithm is an wavefunction-wise
parallelization, based on which the calculation can be done for each
wavefunction separately, and a significant acceleration for the
overall simulation can be expected. However, this is quite nontrivial
for a problem containing eigenvalue problem because of the possible
orthogonalization for all wavefuntions needed during the
simulation. One attractive feature of multilevel correction method is
that the wavefunction-wise parallelization can be partially realized
in the sense that the correction of the wavefunction can be done
individually in the inner iteration. In this work, we have redesigned
the algorithm  to fully take advantage of this feature. It
is worth to rementioning another feature of the multilevel correction
method is that the eigenvalue problem is solved on the corasest mesh,
which can be solved effectively. By combining these two
strategies, a dramatic acceleration for the simulation can be observed
clearly from the numerical experiments.

The outline of this paper is as follows. In Section 2, we recall the
multilevel correction adaptive finite element method for solving
Kohn-Sham equation. In Section 3, we construct an accelerating
multilevel correction adaptive finite element method which can further
improve the solving efficiency for Kohn-Sham equation. In Section 4,
some numerical experiments are presented to demonstrate the efficiency
of the presented algorithm.  Finally, some concluding remarks are
presented in the last section.

\section{Multilevel correction adaptive finite element method for Kohn-Sham equation}
In this section, we review the multilevel correction adaptive finite element method for Kohn-Sham equation (see \cite{HuXieXu}).
To describe the algorithm, we introduce some notation first. Following
\cite{Adams}, we use $W^{s,p}(\Omega)$ to denote Sobolev spaces, and
$\|\cdot\|_{s,p,\Omega}$ and $|\cdot|_{s,p,\Omega}$ to denote the
associated norms and seminorms, respectively. In case $p = 2$, we
denote $H^s(\Omega) = W^{s,2}(\Omega)$ and $H_0^1(\Omega) = \{v \in
H^1(\Omega):v|_{\partial \Omega} = 0\}$, where $v|_{\partial \Omega} =
0$ is in the sense of trace, and denote
$\|\cdot\|_{s,\Omega}=\|\cdot\|_{s,2,\Omega}$. In this paper, we set $V = H_0^1(\Omega)$ for simplicity.

Let $\cH = (H_0^1(\Omega))^N$ be the Hilbert space with the inner product
\begin{eqnarray}
(\Phi,\Psi) = \sum_{i=1}^N\int_{\Omega}\phi_i\psi_idx,
\quad \forall\ \Phi = (\phi_1,\dots,\phi_N),
\ \ \Psi = (\psi_1,\cdots,\psi_N)\in \cH,
\end{eqnarray}
where $\Omega\subset \mathcal{R}^3$ in this paper. For any
$\Psi\in\cH$ and a subdomain $\omega \subset \Omega$, we define
$\rho_{\Psi}=\sum_i|\psi_i|^2$ and
\begin{eqnarray*}
\|\Psi\|_{s,\omega} &=& \left( \sum_{i=1}^{N}
\|\psi_i\|_{s,\omega}^2\right)^{1/2}, \ \ s=0,1.
\end{eqnarray*}
Let $\cQ$ be a subspace of $\cH$ with orthonormality constraints:
\begin{eqnarray}
\cQ = \Big\{  \Psi\in \cH: \Psi^T\Psi = I^{N\times N}    \Big\},
\end{eqnarray}
where $\Phi^T\Psi = \big(\int_{\Omega}\phi_i\psi_jdx\big)_{i,j=1}^N\in \cR^{N\times N}$.

We consider a molecular system consisting of $M$ nuclei with charges
$\{Z_1$, $\cdots$, $Z_M\}$ and locations $\{R_1$, $\cdots$, $R_M\}$,
respectively, and $N$ electrons in the non-relativistic and
spin-unpolarized setting. The general form of Kohn-Sham
energy functional can be demonstrated as follows
\begin{eqnarray}
 E(\Psi) = \int_{\Omega} \left( \frac{1}{2}\sum_{i=1}^N|\nabla \psi_i |^2+V_{\rm ext}(x)\rho_{\Psi}
 +e_{xc}(\rho_{\Psi})\right)dx
 +\frac{1}{2}D(\rho_{\Psi},\rho_{\Psi}),
\end{eqnarray}
for $\Psi = (\psi_1,\psi_2,\cdots,\psi_N)\in \cH$. Here, $V_{\rm ext}$ is
the Coulomb potential defined by $V_{\rm ext} =
-\sum_{k=1}^{M}\frac{Z_k}{|x-R_k|}$, $D(\rho_{\Phi},\rho_{\Phi})$ is
the electron-electron Coulomb energy (Hartree potential) defined by
\begin{eqnarray}
 D(f,g) = \int_{\Omega}f(g*r^{-1})dx = \int_{\Omega}\int_{\Omega}f(x)g(y)\frac{1}{|x-y|}dxdy,
\end{eqnarray}
and $e_{xc}(t)$ is some real function over $[0,\infty)$ denoting the exchange-correlation energy.

The ground state of the system is obtained by solving the minimization problem
\begin{eqnarray}\label{minimization}
 \inf\big\{   E(\Psi):  \Psi \in \cQ     \big\},
\end{eqnarray}
and we refer to \cite{Chen Yang, Cances chakir} for the existence of a minimizer under some conditions.

The Euler-Lagrange equation corresponding to the minimization problem (\ref{minimization}) is the
well known Kohn-Sham equation: Find $(\Lambda, \Phi)\in \cR^{N}\times \cH$ such that
\begin{equation}\label{Kohn_Sham_Equation}
\left\{
\begin{array}{rcl}
H_{\Phi}\phi_i &=&\lambda_i\phi_i\  \ {\rm in}\ \Omega, \quad i=1,\cdots,N,\\
&&\\
\displaystyle\int_\Omega \phi_i \phi_j dx &=&\delta_{ij},
\end{array}
\right.
\end{equation}
where $H_{\Phi}$ is the Kohn-Sham Hamiltonian operator as
\begin{eqnarray}
 H_{\Phi} = -\frac{1}{2}\Delta +V_{\rm ext}
 +\int_{\Omega}\frac{\rho_{\Phi}(y)}{|\cdotp-y|}dy +e'_{xc}(\rho_{\phi})
\end{eqnarray}
with $\Lambda = (\lambda_1,\cdots,\lambda_N)$ and $\lambda_i = (H_\Phi \phi_i,\phi_i)$.
The variational form of the Kohn-Sham equation can be
described as follows: Find $(\Lambda, \Phi)\in \cR^{N}\times \cH$ such that
\begin{equation}\label{Eigenvalue_Problem}
\left\{
\begin{array}{rcl}
(H_{\Phi}\phi_i,v) &=&\lambda_i(\phi_i,v),\quad  \forall v\in H_0^1(\Omega),\quad i=1,2,\cdots,N,\\
  &&\\
\displaystyle\int_\Omega \phi_i \phi_j dx &=&\delta_{ij}.
\end{array}
\right.
\end{equation}
In the ground state of the electronic system, electrons will occupy
the orbitals with the lowest energies. Hence, it corresponds to
finding the left most eigenpairs of Kohn-Sham equation. If the spin
polarization is considered, the number of the eigenpairs is the same
to the number of the electrons. Otherwise, only half eigenpairs are
needed. For simplicity, we only consider the spin-unpolarized case in this paper.

In order to define an efficient way to treat the nonlinear Hartree potential term, we introduce the mixed formulation of the Kohn-Sham equation.
It has been observed from the numerical practice that the term $(V_{\rm Har}\phi_i,v)$ provides the main nonlinearity,
where $V_{\rm Har}$ is the Hartree (electrostatic) potential and its analytical form can be obtained by solving the following Poisson equation:
\begin{eqnarray}
-\Delta V_{\rm Har} &=& 4\pi \rho,
\end{eqnarray}
with a proper Dirichlet boundary condition.

Now, let us define the finite element discretization of
(\ref{Mixed_Eigenvalue_Problem}).  First we generate a shape
regular decomposition $\mathcal{T}_h$ of the computing domain $\Omega$
and let $\mathcal{E}_h$ denote the interior edge set of
$\mathcal{T}_h$. The diameter of a cell $T \in \cT_h$ is denoted by
$h_T$ and the mesh diameter $h$ describes the maximum diameter of all
cells $T \in \cT_h$.  Based on the mesh $\cT_h$, we construct the
linear finite element space denoted by $V_h \subset H_0^1(\Omega)$.
Define $\mathcal H_h = (V_h)^N$ and it is obvious that $\mathcal H_h \subset \cH$.

Then the discrete form of (\ref{Eigenvalue_Problem}) can be described as follows:
Find $(\bar\Lambda_h, \bar\Phi_h)\in \cR^{N}\times \mathcal H_h$ such that
\begin{equation}\label{Nonlinear_Eigenvalue_Problem2 fem}
\left\{
\begin{array}{rcl}
(H_{\bar\Phi_h}\bar\phi_{i,h},v) &=&\bar\lambda_{i,h}(\bar\phi_{i,h},v),
\quad\forall v\in V_h, \quad i=1,\cdots,N,\\
&&\\
\displaystyle\int_\Omega \bar\phi_{i,h} \bar\phi_{j,h} dx &=&\delta_{ij},
\end{array}
\right.
\end{equation}
with $\bar\lambda_{i,h} = (H_{\bar\Phi_h} \bar\phi_{i,h},\bar\phi_{i,h})$.

In order to recover the singularity of Kohn-Sham equation,
the adaptive finite element method (AFEM) is the standard way.
With the adaptive mesh refinement guided by the a posteriori error estimators,
the AFEM can produce an efficient discretization scheme for the singular problems.
The total amount of the mesh elements should be controlled well to make the simulation
continuable and efficient for the Kohn-Sham equation. Based on the above discussion, adaptive
mesh method is a competitive candidate for the refinement strategy.
A standard AFEM process can be described by the following way
\begin{center}
$\cdots$\bf{Solve} $\rightarrow$ \bf{Estimate}
$\rightarrow$ \bf{Mark} $\rightarrow$ \bf{Refine}$\cdots$.
\end{center}
More precisely, to get $\cT_{h_{k+1}}$ from $\cT_{h_k}$, we first solve the discrete equation on $\cT_{h_k}$ to
get the approximate solution and then calculate the \emph{a posteriori} error estimator on each mesh element.
Next we mark the elements with big errors and these elements are refined in such a way that the
triangulation is still shape regular and conforming.

In our simulation, the residual type \emph{a posteriori}
  error estimation is employed to generate the error indicator. First,
we construct the element residual ${\cR}_{T}(\Lambda_h,\Phi_{h})$ and the jump
residual ${\cJ}_e(\Phi_h)$ for the eigenpair approximation
$(\Lambda_h,\Phi_h)$ as follows:
\begin{eqnarray}
&&{\cR}_{T}(\Lambda_h,\Phi_{h}) := \big(H_{\Phi_{h}}\phi_{i,h}-\lambda_{i,h}\phi_{i,h}\big)_{i=1}^{N},
\qquad \text{in } T\in \cT_{h_k}, \\
&&{\cJ}_{e}(\Phi_{h}) := \Big( \frac{1}{2}\nabla \phi_{i,h}|_{T^+}\cdot \nu^+
+\frac{1}{2}\nabla \phi_{i,h}|_{T^-}\cdot \nu^- \Big)_{i=1}^{N},  \quad\text{on } e\in \cE_h,
\end{eqnarray}
where $e$ is the common side of elements $T^+$ and $T^-$ with the unit outward normals $\nu^+$ and
$\nu^-$, respectively. Let $\cE_{h}$ be the set of interior faces (edges or sides) of $\cT_{h}$,
and $\omega_T$ be the union of element sharing a side with $T$.
For $T\in\cT_{h}$, we define the local error indicator ${\eta}_k^2(\Lambda_h,\Phi_{h},T)$ by
\begin{eqnarray}
{\eta}_k^2(\Lambda_h,\Phi_{h},T):=h_T^2\| {\cR}_{T}(\Lambda_h,\Phi_{h})\|_{0,T}^2+\sum_{e\in\cE_{h_k},e
\subset \partial T}h_e\| {\cJ}_{e}(\Phi_{h})\|_{0,e}^2.
\end{eqnarray}
Given a subset $\omega \subset\Omega$, we define the error estimate ${\eta}_k^2(\Lambda_h,\Phi_{h},\omega)$ by
\begin{eqnarray}\label{def of eta and osc}
{\eta}_k^2(\Lambda_h,\Phi_{h},\omega)=\sum_{T\in \cT_{h},T\subset \omega}{\eta}_k^2(\Lambda_h,\Phi_{h},T).
\end{eqnarray}
Based on the error indicator (\ref{def of eta and osc}), we use the
D\"{o}rfler's marking strategy \cite{Dorfler} to mark all elements in $\cM_k$ for local refinement.

Since solving large-scale nonlinear eigenvalue problem is quite time-consuming compared to that of boundary value problem, a multilevel correction adaptive method for solving Kohn-Sham equation
was designed in \cite[Algorithms 1 and 3]{HuXieXu}.
The multilevel correction method transforms solving Kohn-Sham equation
on the adaptive refined mesh $\mathcal T_h$ to the solution
$\widetilde \Phi_h$ of the associated linear boundary value problems
on $\mathcal T_h$ and nonlinear eigenvalue problem in a fixed low
dimensional subspace $V_{H,h}=V_H+{\rm span}\{\widetilde
  \Phi_h\}$ which is build with the finite element space $V_H$ on the
coarse mesh $\mathcal T_H$ and $N$ finite element functions
$\widetilde\Phi_h = (\widetilde\phi_{i,h}, \cdots,
\widetilde\phi_{N,h})$.  In the multilevel correction adaptive method
for solving Kohn-Sham equation, we need to solve $N$ linear boundary
value problems in each adaptive space and a small-scale Kohn-Sham
equation in the correction space $V_{H,h}$.  Since there is no
eigenvalue problem solving in the fine mesh $V_h$, the multilevel
correction adaptive method has a better efficiency than the direct
AFEM.  The dimension of the correction space $V_{H,h}$ is fixed and
small in the multilevel correction adaptive finite element method.
But we need to solve a nonlinear eigenvalue problem in the correction
space $V_{H,h}$. Always, some type of nonlinear iteration steps are
required to solve this nonlinear eigenvalue problems. When the system
includes large number of electrons, the number of required nonlinear
iteration steps is always very large.  Furthermore, the correction
space contains $N$ basis functions $\widetilde\Phi_h$ defined in the
fine space $V_h$.  In order to guarantee the calculation accuracy, we
need to assemble matrices in the fine space when it comes to the basis
functions of ${\rm span}\{\widetilde{\Phi}_{h}\}$ defined on the fine
mesh $\mathcal T_h$.  This part of work depends on the number of
electrons and dimension of $V_h$ as $N^2\times {\rm dim}V_h$.  In
order to improve the efficiency of the correction step, a new type of
nonlinear iteration method and an eigenpairwise parallel correction method with
efficient implementing techniques will be designed in the next
section.

\section{Further improvement of the method towards the efficiency}
In this section, we give a new type of AFEM
which is a combination of the multilevel correction scheme, nonlinearity separating technique and an efficient parallel
method for the Hartree potential. The nonlinearity separating technique is designed based on the different performance of the exchange-correlation and Hartree potentials.
We will use the self consistent field (SCF) iteration steps to treat the weak nonlinear exchange-correlation potential and
the multilevel correction method for the strong nonlinear Hartree potential. Furthermore, combining the eigenpairwise parallel idea from \cite{XuXieZhang} and
the special structure of the Hartree potential, the Hartree potential can be treated by the parallel way and an efficient implementing technique.

\subsection{A strategy of separating nonlinear terms}
The nonlinearity of the Kohn-Sham equation comes from Hartree and exchange-correlation potentials.
As discussed in Section 1, the Hartree potential has a strong nonlinearity which leads the main nonlinearity
of the Kohn-Sham equation,
while the nonlinearity of the exchange-correlation potential is weak.
Based on such a property, we modify the standard SCF iteration into a nested iteration which includes outer SCF iteration
for the exchange-correlation potential, and inner multilevel correction iteration for the Hartree potential.
Furthermore, we will also design an eigenpairwise parallel multilevel correction method for solving the nonlinear eigenvalue problems
associated with the Hartree potential in the inner iterations.


Actually, this section is to define an efficient numerical method to
solve the Kohn-Sham equation on the refined mesh $\cT_{h_{k+1}}$ with the proposed
nested iteration here.
Different from the standard SCF iteration, we decompose the SCF iteration into outer iteration for exchange-correlation potential
and inner iteration for the Hartree potential.
This type of nonlinear treatment is based on the understanding of the nonlinearity strengthes from the
Hartree and exchange-correlation potentials, respectively. In the practical models and computing experience,
the nonlinearity of Hartree Potential is stronger than that of the exchange-correlation potential.

Given an initial value $(\Lambda_{h_{k+1}}^{(\ell)}, \Phi_{h_{k+1}}^{(\ell)})$ for the Kohn-Sham equation, we should solve the following nonlinear eigenvalue
problem in each step of the outer SCF iteration method for the exchange-correlation potential:
For $i=1,\cdots,N$,
find $(\lambda_{i,h_{k+1}}^{(\ell+1)},\phi_{i,h_{k+1}}^{(\ell+1)})\in \mathcal R\times V_{h_{k+1}}$
such that
\begin{eqnarray}\label{Kohn_Sham_EC}
\begin{array}{r}
L\left(\Phi_{h_{k+1}}^{(\ell)};\phi_{i,h_{k+1}}^{(\ell+1)}, \gamma_{h_{k+1}}\right)
= \lambda_{i,h_{k+1}}^{(\ell+1)}\left(\phi_{i,h_{k+1}}^{(\ell+1)},\gamma_{h_{k+1}}\right),\ \ \
 \forall \gamma_{h_{k+1}}\in V_{h_{k+1}},
\end{array}
\end{eqnarray}
where
\begin{eqnarray}
\begin{array}{r}
L\left(\Phi_{h_{k+1}}^{(\ell)};\phi_{i,h_{k+1}}^{(\ell+1)}, \gamma_{h_{k+1}}\right)=a(\phi_{i,h_{k+1}}^{(\ell+1)}, \gamma_{h_{k+1}})+D\left(\rho_{\Phi_{h_{k+1}}^{(\ell+1)}},\phi_{i,h_{k+1}}^{(\ell+1)}\gamma_{h_{k+1}}\right)
+\left(V_{xc}(\rho_{\Phi_{h_{k+1}}^{(\ell)}})\phi_{i,h_{k+1}}^{(\ell+1)},\gamma_{h_{k+1}}\right),
\end{array}
\end{eqnarray}
and
\begin{eqnarray}\label{Bilinear_Form_a}
a(\phi_{i,h_{k+1}}^{(\ell+1)}, \gamma_{h_{k+1}}) = \frac{1}{2}\left(\nabla \phi_{i,h_{k+1}}^{(\ell+1)}, \nabla \gamma_{h_{k+1}}\right) + \left(V_{\rm ext}\phi_{i,h_{k+1}}^{(\ell+1)},\gamma_{h_{k+1}}\right).
\end{eqnarray}

In the above SCF iteration, only the exchange-correlation potential is linearized.
Thus, we still need to solve a nonlinear eigenvalue problem in each iteration step, and the nonlinearity is caused by the Hartree potential.
Because the nonlinearity of the exchange-correlation potential is weak, so only a few steps are needed for the above SCF iteration.
Further, for the involved nonlinear eigenvalue problem whose nonlinearity is caused by the Hartree potential, we can design a
new strategy to solve it efficiently.

It is obvious that equation (\ref{Kohn_Sham_EC}) is still a nonlinear eigenvalue problem with the nonlinear term of Hartree potential.
Now, we spend the main attentions to design an efficient eigenpairwise parallel way to solve the nonlinear eigenvalue problem (\ref{Kohn_Sham_EC}).
The idea and method here come from \cite{HuXieXu,XuXieZhang,ZhangXuXie}. This type of method is built based
on the low dimensional space defined on the coarse mesh $\cT_H$. But the method in this paper is the eigenpairwise
parallel way to implement the augmented subspace method which is different from \cite{HuXieXu}.
Compared with the standard SCF iteration method for Kohn-Sham equation, there is no inner products for orthogonalization process in the high dimensional space $V_{h_{k+1}}$, which is always the bottle neck for parallel computing.

In order to define the eigenpairwise multilevel correction method for the Hartree potential, we transform the Kohn-Sham equation (\ref{Eigenvalue_Problem})
into the following equivalently mixed form: 
Find $(\Lambda, \Phi, w)\in \mathcal R^{N}\times  \cH \times V_\Gamma$ such that
\begin{eqnarray}\label{Mixed_Eigenvalue_Problem}
\left\{
\begin{array}{rcl}
a(\phi_i,\gamma)+(w\phi_i,\gamma)+(V_{xc}\phi_i,\gamma)
&=&\lambda_i(\phi_i,\gamma),\quad \ \  \forall \gamma\in V,\\
(\nabla w, \nabla v) &=& 4\pi(\rho_\phi, v),\quad\quad \ \forall v\in V,\\
\displaystyle\int_\Omega \phi_i \phi_j dx &=&\delta_{ij},
\end{array}
\right.
\end{eqnarray}
where the function set $V_\Gamma$ is defined by the trace of Hartree potential on the boundary $\partial\Omega$
\begin{eqnarray*}
V_\Gamma = \left\{v\in H^1(\Omega)\ \Big| v_{\partial\Omega}=\int_{\mathcal R^3}\frac{\rho_\Phi(y)}{|x-y|}dy\Big|_{\partial\Omega}\right\}.
\end{eqnarray*}

In the practical computation, the discrete form of (\ref{Mixed_Eigenvalue_Problem}) can be described as follows:
Find $(\bar\Lambda_h, \bar\Phi_h, \bar w_h)\in \cR^{N}\times \cH_h \times V_{\Gamma, h}$ such that
\begin{eqnarray}\label{Mixed_Eigenvalue_Problem_Discrete}
\left\{
\begin{array}{rcl}
a( \bar\phi_{i,h}, \gamma_h)+(\bar w_h\bar \phi_{i,h},\gamma_h)-(V_{xc}\bar\phi_{i,h},\gamma_h)
&=&\bar\lambda_{i,h}(\bar\phi_{i,h},\gamma_h),\quad  \forall \gamma_h\in V_h,\\
(\nabla \bar w_h, \nabla v_h)- 4\pi \sum\limits_{i=1}^N(\bar\phi_{i,h}^2, v_h)&=&0,\quad\quad\quad \ \ \forall v_h\in V_h,\\
\displaystyle\int_\Omega \bar\phi_{i,h} \bar\phi_{j,h} dx &=&\delta_{ij},
\end{array}
\right.
\end{eqnarray}
where the set $V_{\Gamma,h}$ includes linear finite element functions with the boundary condition which is computed by the fast multipole method \cite[(29)]{BaoHuLiu}.

The corresponding scheme is defined by Algorithm \ref{Parallel_Aug_Subspace_Method}, where we can
find that the augmented subspace method transforms solving the nonlinear eigenvalue problem
into the solution of the boundary value problem and some small scale nonlinear eigenvalue problems in the low dimensional space $S_{H,h_{k+1}}$.

\begin{algorithm}[htbp]
\caption{An accelerating multilevel correction adaptive finite element method}\label{Parallel_Aug_Subspace_Method}
\DontPrintSemicolon
\BlankLine
Solve the Kohn-Sham equation in the initial finite element space $V_{h_1}$:
Find $(\Lambda_{h_1},\Phi_{h_1})\in \cR^{N} \times V_{h_1}$ such that
\begin{eqnarray*}
(\cH_{\Phi_{h_1}}\phi_{i,h_1},v_{h_1})=\lambda_{i,h_1}(\phi_{i,h_1},v_{h_1}),
\quad \forall v_{h_1} \in V_{h_1}, \quad i=1,2,\cdots,N.
\end{eqnarray*}\;
Set $k=1$.\;
Compute the local error indicators ${\eta}_k(\Phi_{h_k},T)$ for each element $T\in\mathcal T_{h_k}$,
and construct a new mesh $\cT_{h_{k+1}}$ according to ${\eta}_k(\Phi_{h_k},T)$ and the Dorfler's Marking Strategy.\;
Set $\ell=0$ and the initial value $\Lambda_{h_{k+1}}^{(\ell)} = \Lambda_{h_k}$,
$\Phi_{h_{k+1}}^{(\ell)}=\Phi_{h_k}$, $w_{h_{k+1}}^{(\ell)}=w_{h_k}$.\;
Define the following linear boundary value problem:
Find $\widehat{\phi}_{i,h_{k+1}}^{(\ell+1)} \in V_{h_{k+1}}, i=1,\cdots,N$
such that
\begin{eqnarray}\label{Aux_Linear_Problem}
\begin{array}{rcr}
&&a(\widehat{\phi}_{i,h_{k+1}}^{(\ell+1)},\gamma_{h_{k+1}})
+(\widehat{w}_{h_{k+1}}^{(\ell)}\widehat\phi_{i,h_{k+1}}^{(\ell+1)},\gamma_{h_{k+1}})+
\left(V_{xc}(\rho_{\Phi_{h_{k+1}}^{(\ell)}})\widehat\phi_{i,h_{k+1}}^{(\ell+1)},\gamma_{h_{k+1}}\right)
 =\lambda_{i,h_{k+1}}^{(\ell)}(\phi_{i,h_{k+1}}^{(\ell)},\gamma_{h_{k+1}}).
\end{array}
\end{eqnarray}
Solve this elliptic equation with some type of iteration method to obtain an approximation $\widetilde{\phi}_{i,h_{k+1}}^{(\ell+1)}$.\;
Define the following elliptic problem: Find $\widehat w_{h_{k+1}}^{(\ell+1)} \in V_{h_{k+1}}$ such that
\begin{eqnarray}\label{W_h_equation}
(\nabla \widehat w_{h_{k+1}}^{(\ell+1)}, \nabla v_{h_{k+1}})= 4\pi \sum_{i=1}^N\left((\widetilde\phi_{i,h_{k}}^{(\ell+1)})^2, v_{h_{k+1}}\right),
\ \ \ \forall v_{h_{k+1}}\in  V_{h_{k+1}}.
\end{eqnarray}
Solve this elliptic equation with some type of iteration method to obtain an approximation
$\widetilde{w}_{h_{k+1}}.$\;
For $i=1,\cdots,N$, define  $S_{H,h_{k+1}}=V_H+{\rm span}\{\widetilde{\phi}_{i,h_{k+1}}^{(\ell+1)}\}$ and
$W_{H,h_{k+1}} = V_H +{\rm span}\{\widetilde w_{h_{k+1}}^{(\ell+1)}\}$.
Solve the following nonlinear eigenvalue problem:
Find $(\lambda_{i,h_{k+1}}^{(\ell+1)},\phi_{i,h_{k+1}}^{(\ell+1)}, w_{h_{k+1}}^{(\ell+1)}) \in \mathcal R\times S_{H,h_{k+1}}\times W_{H,h_{k+1}}$
such that
\begin{equation}\label{Nonlinear_Eig_Hh_Full}
\left\{
\begin{array}{r}
L(\Phi_{h_{k+1}}^{(\ell)};\phi_{i,h_{k+1}}^{(\ell+1)}, \gamma_{H,h_{k+1}})
= \lambda_{i,h_{k+1}}^{(\ell+1)}(\phi_{i,h_{k+1}}^{(\ell+1)},\gamma_{H,h_{k+1}}), \ \ \ \forall \gamma_{H,h_{k+1}} \in S_{H,h_{k+1}}\\
(\nabla w_{h_{k+1}}^{(\ell+1)}, \nabla v_{H,h_{k+1}}) =4\pi \sum\limits_{i=1}^N\left((\phi_{i,h_{k+1}}^{(\ell+1)})^2, v_{H,h_{k+1}}\right), \ \ \ \forall v_{H,h_{k+1}} \in W_{H,h_{k+1}}.
\end{array}
\right.
\end{equation}
\;
If $\|\rho_{\Phi_{h_{k+1}}^{(\ell+1)}}-\rho_{\Phi_{h_{k+1}}^{(\ell)}}\|_1\geq {\rm tol}$, set $\ell=\ell+1$ and go to step 5, else go to step 9.\;
Set $\left(\Lambda_{h_{k+1}},\Phi_{h_{k+1}},w_{h_{k+1}}\right)
=\left(\Lambda_{h_{k}}^{(\ell+1)},\Phi_{h_{k+1}}^{(\ell+1)},w_{h_{k+1}}^{(\ell+1)}\right)$, $k=k+1$ and go to step 3.
\end{algorithm}

The linear boundary value problems (\ref{Aux_Linear_Problem})  and (\ref{W_h_equation})
in Algorithm \ref{Parallel_Aug_Subspace_Method}
are easy to be solved by the efficient and matured linear solvers.
For example, the linear solvers based on the multigrid method always give the optimal efficiency
for solving the second order elliptic problems.
Compared with the multilevel correction method defined in \cite{HuXieXu}, the
obvious difference is that Algorithm \ref{Parallel_Aug_Subspace_Method} treats the different orbit
associated with $\widetilde\phi_{i,h_{k+1}}$, $i=1, \cdots, N$, indepdently.
This way can avoid doing inner products for orthogonalization  in the high dimensional space $V_{h_{k+1}}$.
Furthermore, we will design an efficient numerical method for solving the nonlinear eigenvalue problem (\ref{Nonlinear_Eig_Hh_Full})
with a special implemenmting techqniue in the next subsection.

\subsection{New algorithm and its parallel implementation}
In this subsection, we further study two important characteristics of Algorithm \ref{Parallel_Aug_Subspace_Method}, which will help to
accelerate the solving efficiency. The first characteristic is that the SCF iteration for the nonlinear eigenvalue
problem (\ref{Nonlinear_Eig_Hh_Full}) can be implemented efficiently. The idea here comes from \cite{XuXieZhang,ZhangXuXie}.
The second characteristic is that Algorithm \ref{Parallel_Aug_Subspace_Method} is suitable for eigenpairwise parallel computing.

We begin with the first characteristic of Algorithm \ref{Parallel_Aug_Subspace_Method}
by introducing an efficient implementation of the iteration method for
the nonlinear eigenvalue problem (\ref{Nonlinear_Eig_Hh_Full}). The corresponding implementation scheme
is described by Algorithm \ref{SCF_Iteration_Hartree}.
\begin{algorithm}[htbp]
\caption{SCF iteration for step 7 in Algorithm \ref{Parallel_Aug_Subspace_Method}}\label{SCF_Iteration_Hartree}
\DontPrintSemicolon
\BlankLine
Set $s=1$ and the initial value $\Phi_{h_{k+1}}^{(\ell+1,s)}=\widetilde\Phi_{h_k}^{(\ell+1)}$, $w_{h_{k+1}}^{(\ell+1,s)}=\widetilde w_{h_k}^{(\ell+1)}$.
Do the following nonlinear iteration:
\begin{enumerate}
\item[(a).]  For $i=1,2,\cdots,N$, solve the following linear eigenvalue problem.
\begin{eqnarray}\label{Eigenvalue_Hh_Fullalg2}
&&a(\phi_{i,h_{k+1}}^{(\ell+1,s+1)},\gamma_{H,h_{k+1}})
+(w_{h_{k+1}}^{(\ell+1,s)}\phi_{i,h_{k+1}}^{(\ell+1,s+1)},\gamma_{H,h_{k+1}})
+\left(V_{xc}(\rho_{\Phi_{h_{k+1}}^{(\ell)}})\phi_{i,h_{k+1}}^{(\ell+1,s+1)},\gamma_{H,h_{k+1}}\right)\nonumber\\
&&=\lambda_{i,h_{k+1}}^{(\ell+1,s+1)}(\phi_{i,h_{k+1}}^{(\ell+1,s+1)},\gamma_{H,h_{k+1}}),\ \ \ \forall \gamma_{H,h_{k+1}}\in S_{H,h_{k+1}}.
\end{eqnarray}
Solve (\ref{Eigenvalue_Hh_Fullalg2}) and choose the eigenfunction $\phi_{i,h_{k+1}}^{(\ell+1,s+1)}$ that
has the largest component in span$\{\widetilde{\phi}_{i,h_{k+1}}^{(\ell+1)}\}$ among all the eigenfunctions.
\item [(b).] Solve the following elliptic problem: Find $w_{h_{k+1}}^{(\ell+1,s+1)} \in W_{H,h_{k+1}}$ such that
\begin{eqnarray}\label{hartree}
\hskip-2cm
\big(\nabla w_{h_{k+1}}^{(\ell+1,s+1)}, \nabla v_{H,h_{k+1}}\big) = 4\pi \sum_{i=1}^N\big((\phi_{i,h_{k+1}}^{(\ell+1,s+1)})^2,v_{H,h_{k+1}}\big),
\  \forall v_{H,h_{k+1}}\in W_{H,h_{k+1}}.
\end{eqnarray}
\item [(c).] If $\|\Phi_{h_{k+1}}^{(\ell+1,s+1)}-\Phi_{h_{k+1}}^{(\ell+1,s)}\|_1 < {\rm tol}$,
stop the iteration. Else set $s=s +1$ and go to (a) to do the next iteration step.
\end{enumerate}\;
We define the output as the final eigenpair approximation 
$$\left(\Lambda_{h_{k+1}}^{(\ell+1)},\Phi_{h_{k+1}}^{(\ell+1)},w_{h_{k+1}}^{(\ell+1)}\right)
=\left(\Lambda_{h_{k+1}}^{(\ell+1,s+1)},\Phi_{h_{k+1}}^{(\ell+1,s+1)},w_{h_{k+1}}^{(\ell+1,s+1)}\right).$$
\end{algorithm}

Different from the standard SCF method, the one in Algorithm \ref{SCF_Iteration_Hartree}
use the mixed form of Kohn-Sham equation to do the nonlinear iteration.
Even the mixed form is equivalent to the standard one, it provides a chance to design an efficient
implementing way to do the nonlinear iteration.
This means the remaining part of this subsection is to discuss how to perform Algorithm \ref{SCF_Iteration_Hartree} efficiently based on the special structure of the correction
spaces $S_{H,h}$ and $W_{H,h}$. The designing process for the efficient implementing way
also shows a reason to treat the Hartree potential and exchange-correlation potential separately.

From the definition of Algorithm \ref{SCF_Iteration_Hartree}, we can find that solving eigenvalue problem
(\ref{Eigenvalue_Hh_Fullalg2}) and linear boundary value problem (\ref{hartree}) need very small computation work
since the dimensions of $S_{H,h_{k+1}}$ and $W_{H,h_{k+1}}$ are very small.
But both $S_{H,h_{k+1}}$ and $W_{H,h_{k+1}}$ include finite element functions defined on the finer
mesh $\mathcal T_{h_{k+1}}$. In order to guarantee the accuracy, we need to assemble the
matrices and right hand side terms in  (\ref{Eigenvalue_Hh_Fullalg2}) and (\ref{hartree})
in the finer mesh $\mathcal T_{h_{k+1}}$ which needs the computational work $\mathcal O({\rm dim}V_{h_{k+1}})$.
Based on this consideration, the key point for implementing Algorithm \ref{SCF_Iteration_Hartree}
efficiently is to design an efficient way to assemble the concerned matrices and right hand side terms.
For the description of implementing technique here, let  $\{\psi_{k,H}\}_{1\leq k\leq N_H}$ denotes the Lagrange basis function for the coarse finite element space $V_H$.

%
In order to show the main idea here, let us consider the matrix version of eigenvalue problem (\ref{Eigenvalue_Hh_Fullalg2}) as follows
\begin{equation}\label{Eigenvalue_H_h}
\left(
\begin{array}{cc}
A_H & b_{Hh}\\
b_{Hh}^T&\beta
\end{array}
\right)
\left(
\begin{array}{c}
\boldsymbol{\upphi}_{i,H} \\
\theta_{2,i}
\end{array}
\right)
=\lambda_i\left(
\begin{array}{cc}
M_H & c_{Hh}\\
c_{Hh}^T&\gamma
\end{array}
\right)
\left(
\begin{array}{c}
\boldsymbol{\upphi}_{i,H}\\
\theta_{2,i}
\end{array}
\right),
\end{equation}
where $A_H \in \mathbb R^{N_H\times N_H}$, $b_{Hh} \in \mathbb R^{N_H}$, $\boldsymbol{\upphi}_{i,H}\in \mathbb R^{N_H}$ and $\beta$, $\gamma$, $\theta_{2,i}\in\mathbb R$.

Since it is required to solve $N$ linear eigenvalue problems which can be assembled in the same way,
we use $\widetilde w_h$ and $\widetilde \phi_{h}$ to denote $\widetilde w_{h_{k+1}}$ and $\widetilde \phi_{i,h_{k+1}}$, respectively, for simplicity.
We should know that the following process is performed for each SCF iteration.
Here, for readability, the upper indices used in our algorithm is also omitted.

It is obvious that the matrix $M_H$, the vector $c_{Hh}$ and the scalar $\gamma$ will not change
during the nonlinear iteration process as long as we have obtained the function $\widetilde \phi_h$.
But the matrix $A_H$, the vector $b_{Hh}$ and the scalar $\beta$ will change during the nonlinear
iteration process. Then it is required to consider the efficient implementation to update the
the matrix $A_H$, the vector $b_{Hh}$ and the scalar $\beta$ since there is a function $\widetilde \phi_{h}$
which is defined on the fine mesh $\mathcal T_h$. The aim here is to propose an efficient method to update
the matrix $A_H$, the vector $b_{Hh}$ and the scalar $\beta$ without computing on the fine mesh $\mathcal T_h$
during the nonlinear iteration process.
Assume we have a given initial value $V_{\rm Har}=w_H+\theta_1\widetilde w_h$ for Hartree potential. Now, in order to carry out the nonlinear iteration
for the eigenvalue problem (\ref{Eigenvalue_H_h}), we come to consider the computation for the matrix $A_H$,
the vector $b_{Hh}$ and the scalar $\beta$.

From the definitions of the space $S_{H,h}$ and the eigenvalue problem (\ref{Eigenvalue_Hh_Fullalg2}), the matrix
$A_H$ has the following expansion
\begin{eqnarray}\label{Adef}
(A_H)_{j,k}&=&\frac{1}{2}\int_\Omega\nabla\psi_{k,H}\nabla \psi_{j,H}dx +\int_\Omega V_{\rm ext}\psi_{k,H}\psi_{j,H}dx +\int_\Omega V_{\rm Har}\psi_{k,H}\psi_{j,H}dx +\int_{\Omega}V_{xc}\psi_{k,H}\psi_{j,H}dx \nonumber\\
&=&\frac{1}{2}\int_\Omega\nabla\psi_{k,H}\nabla \psi_{j,H}dx +\int_\Omega V_{\rm ext}\psi_{k,H}\psi_{j,H}dx +\int_\Omega (w_H+\theta_1\widetilde w_h)\psi_{k,H}\psi_{j,H}dx +\int_{\Omega}V_{xc}\psi_{k,H}\psi_{j,H}dx \nonumber\\
&=&\Big(\frac{1}{2}\int_\Omega\nabla\psi_{k,H}\nabla \psi_{j,H}dx +\int_\Omega V_{\rm ext}\psi_{k,H}\psi_{j,H}dx \Big)+\int_\Omega w_H\psi_{k,H}\psi_{j,H}dx
+\theta_1\int_\Omega\widetilde w_h\psi_{k,H}\psi_{j,H}dx  \nonumber\\
&&+\int_{\Omega}V_{xc}\psi_{k,H}\psi_{j,H}dx \nonumber\\
&:= &(A_{H,1})_{j,k} + (A_{H,2})_{j,k} + \theta_1(A_{H,3})_{j,k}+ (A_{H,4})_{j,k},\ \ \ 1\leq j,k\leq N_H.
\end{eqnarray}
Here $A_{H,1}$, $A_{H,3}$ remain unchanged during the inner nonlinear iteration.
During the inner loop, the exchange-correlation potential $V_{xc}$ remains unchanged, so the matrix $A_{H,4}$ also remains unchanged during the inner nonlinear iteration.
Thus we can assemble $A_{H,1}$, $A_{H,3}$, $A_{H,4}$ in advance, and call these data directly in each iteration step.
The matrix $A_{H,2}$ will change during the nonlinear iteration because $w_H$ will change.
But $A_{H,2}$ is defined on the coarse space $V_H$, which can be assembled by the small computational work $\mathcal O(N_H)$.

The vector $b_{Hh}$ has the following expansion
\begin{eqnarray}\label{bdef}
(b_{Hh})_{j}&=&\frac{1}{2}\int_\Omega\nabla \widetilde \phi_{h}\nabla\psi_{j,H}dx +\int_\Omega V_{\rm ext}\widetilde \phi_{h}\psi_{j,H}dx
+\int_\Omega V_{\rm Har}\widetilde \phi_{h}\psi_{j,H}dx +\int_{\Omega}V_{xc}\widetilde \phi_{h}\psi_{j,H}dx \nonumber\\
&=&\frac{1}{2}\int_\Omega\nabla \widetilde \phi_{h}\nabla\psi_{j,H}dx +\int_\Omega V_{\rm ext}\widetilde \phi_{h}\psi_{j,H}dx
+\int_\Omega (w_H+\theta_1\widetilde w_h)\widetilde \phi_{h}\psi_{j,H}dx +\int_{\Omega}V_{xc}\widetilde \phi_{h}\psi_{j,H}dx \nonumber\\
&=&\Big(\frac{1}{2}\int_\Omega\nabla \widetilde \phi_{h}\nabla\psi_{j,H}dx +\int_\Omega V_{\rm ext}\widetilde \phi_{h}\psi_{j,H}dx \Big)
+\int_\Omega w_H\widetilde \phi_{h}\psi_{j,H}dx +\theta_1\int_\Omega\widetilde w_h\widetilde \phi_{h}\psi_{j,H}dx \nonumber\\
&&+\int_{\Omega}V_{xc}\widetilde \phi_{h}\psi_{j,H}dx \nonumber\\
&:=& (b_{H,1})_{j} + (b_{H,2})_{j} + \theta_1(b_{H,3})_{j}+ (b_{H,4})_{j},\ \ \ 1\leq j\leq N_H.
\end{eqnarray}
Because $\widetilde \phi_h$ is a basis function of the correction space $S_{H,h}$, so $b_{H,1}, b_{H,3}, b_{H,4}$ remain unchanged during the inner nonlinear iteration.
So we can compute these three vectors in advance.
In order to compute  $b_{H,2}$, we first assemble a matrix $A_{h4}$ in the fine space $V_h$
\begin{eqnarray}\label{Ah4def}
(A_{h4})_{j,k}=(\widetilde\phi_{h}\psi_{k,H}, \psi_{j,H}), \ \ \ \ 1\leq j,k\leq N_H,
\end{eqnarray}
which will be fixed during the inner nonlinear iteration and can be computed in advance.
Let us define $w_H=\sum_{k=1}^{N_H}w_k\psi_{k,H}$, and ${\bf w}_H=(w_1,\cdots,w_{N_H})^T$. Then, $b_{H,2}=A_{h4}{\bf w}_H$.

The scalar $\beta$ has the following expansion
\begin{eqnarray}\label{betadef}
\beta&=&\frac{1}{2}\int_\Omega\nabla\widetilde \phi_{h}\nabla \widetilde \phi_{h}dx +\int_\Omega V_{\rm ext}\widetilde \phi_{h}\widetilde \phi_{h}dx
+\int_\Omega V_{\rm Har}\widetilde \phi_{h}\widetilde \phi_{h}dx +\int_{\Omega}V_{xc}\widetilde \phi_{h}\widetilde \phi_{h}dx \nonumber\\
&=&\frac{1}{2}\int_\Omega\nabla\widetilde \phi_{h}\nabla \widetilde \phi_{h}dx +\int_\Omega V_{\rm ext}\widetilde \phi_{h}\widetilde \phi_{h}dx
+\int_\Omega (w_H+\theta_1\widetilde w_h)\widetilde \phi_{h}\widetilde \phi_{h}dx +\int_{\Omega}V_{xc}\widetilde \phi_{h}\widetilde \phi_{h}dx \nonumber\\
&=&\Big(\frac{1}{2}\int_\Omega\nabla\widetilde \phi_{h}\nabla \widetilde \phi_{h}dx +\int_\Omega V_{\rm ext}\widetilde \phi_{h}\widetilde \phi_{h}dx \Big)
+\int_\Omega w_H\widetilde \phi_{h}\widetilde \phi_{h}dx +\theta_1\int_\Omega\widetilde w_h\widetilde \phi_{h}\widetilde \phi_{h}dx \nonumber\\
&&+\int_{\Omega}V_{xc}\widetilde \phi_{h}\widetilde \phi_{h}dx \nonumber\\
&:=& \beta_1 + \beta_2 + \theta_1\beta_3+ \beta_4.
\end{eqnarray}
Here $\beta_1$, $\beta_3$, $\beta_4$ remain unchanged during the inner nonlinear iteration.
In order to compute  $\beta_2$ efficiently, we first assemble a vector ${\tt rhs}$ in the fine space $V_h$:
\begin{eqnarray}\label{rhsdef}
({\tt rhs})_{j}=(\widetilde\phi_{h}\widetilde\phi_{h}, \psi_{j,H}),\ \ \ 1\leq j\leq N_H.
\end{eqnarray}
It is obvious that the vector ${\tt rhs}$ will be fixed during the inner nonlinear iteration
and $\beta_2={\tt rhs}^T{\bf w}_H$.

Next, we consider the efficient scheme for solving the Hartree potential equation (\ref{hartree}).
Based on the structure of the space $W_{H,h}$, the matrix version of (\ref{hartree}) can be written as follows
\begin{equation}\label{Hartree_H_h}
\left(
\begin{array}{cc}
C_H & d_{Hh}\\
d_{Hh}^T&\zeta
\end{array}
\right)
\left(
\begin{array}{c}
\mathbf{w}_{H} \\
\theta_1
\end{array}
\right)
=4\pi\left(
\begin{array}{c}
f_{H}\\
g
\end{array}
\right),
\end{equation}
where $C_H \in \mathbb R^{N_H\times N_H}$, $d_{Hh} \in \mathbb R^{N_H}$, $f_H\in \mathbb R^{N_H}$, $\zeta\in\mathbb R$ and $g\in\mathbb R$.

It is obvious that the matrix $(C_H)_{j,k}=(\nabla\psi_{k,H},\nabla\psi_{j,H})$, the vector
$(d_{Hh})_j=(\nabla \widetilde w_h, \nabla\psi_{j,H})$ and the
scalar $\gamma=(\nabla\widetilde w_h,\nabla\widetilde w_h)$ will not change
during the nonlinear iteration process as long as we have obtained the function $\widetilde w_h$. So the main task is to assemble the right hand side term $[f_H, g]^T$.
Since the density function is the sum of the square of $N$ approximate eigenfunctions,
we should do the summation according to the lower index $i$ in our description.

For the right hand term $f_H$ of (\ref{Hartree_H_h}), we have
\begin{eqnarray}\label{fH}
(f_{H})_j&=&\left(\sum_{i=1}^N(\phi_{i,H}+\theta_{2,i}\widetilde\phi_{i,h})^2,\psi_{j,H}\right)\nonumber\\
&=&\sum_{i=1}^N\left((\phi_{i,H}^2+2\theta_{2,i}\phi_{i,H}\widetilde\phi_{i,h}
+\widetilde\phi_{i,h}^2),\psi_{j,H}\right)\nonumber\\
&:=& \sum_{i=1}^N\left((f_{H,i,1})_{j}  + \theta_{2,i}(f_{H,i,2})_{j}+ (f_{H,i,3})_{j}\right),\ \ \ 1\leq j\leq N_H.
\end{eqnarray}
The vector $f_{H,i,3}$ remains unchanged during the inner SCF iteration and $f_{H,i,3}$ is just the vector
${\tt rhs}$ in (\ref{rhsdef}) assembled
for the $i$-th linear boundary value problem. Because $\phi_{i,H}$ will change during the inner nonlinear iteration, so $f_{H,i,1}$ needs to be computed during the iteration.
Fortunately, $f_{H,i,1}$ is defined on the coarse space $V_H$, which only needs the computational work $\mathcal O(N_H)$.
The vector $f_{H,i,2}$ can be computed in the way $f_{H,i,2}=A_{h4}\boldsymbol{\upphi}_{i,H}$,
where $\boldsymbol{\upphi}_{i,H}$ is the coefficient vector of $\phi_{i,H}$ with respect to the basis functions of $V_H$.

From the structure of $W_{H,h}$, the scalar $g$ has the following expansion
\begin{eqnarray}\label{gdef}
g&=&\left(\sum_{i=1}^N(\phi_{i,H}+\theta_{2,i}\widetilde\phi_{i,h})^2,\widetilde w_{h}\right)\nonumber\\
&=&\sum_{i=1}^N\left((\phi_{i,H}^2+2\theta_{2,i}\phi_{i,H}\widetilde\phi_{i,h}+\widetilde\phi_{i,h}^2),\widetilde w_{h}\right)\nonumber\\
&:=& \sum_{i=1}^N\big(g_{i,1}  + \theta_{2,i}g_{i,2}+ g_{i,3}\big).
\end{eqnarray}

The scalar $g_{i,3}$ remains unchanged during the inner nonlinear iteration.
In order to assemble $g_{i,1}$ and $g_{i,2}$, let us define a matrix $G\in \mathbb R^{N_H\times N_H}$ and a vector $F_i\in \mathbb R^{N_H}$
in the following way which remain unchanged during the inner nonlinear iteration:
\begin{eqnarray}\label{FG}
(G)_{j,k}=(\widetilde w_{h}\psi_{k,H},\psi_{j,H}) \ \ \
\text{and} \ \ \
(F_i)_j=(\widetilde w_{h}\widetilde \phi_{i,h},\psi_{j,H}), \ \ 1\leq j,k\leq N_H, \ \ i=1,\cdots,N.
\end{eqnarray}
Then
\begin{eqnarray}
g_{i,1}=\boldsymbol{\upphi}_{i,H}^TG \boldsymbol{\upphi}_{i,H} \ \ \  \text{and} \ \ \ g_{i,2}=F_i^T\boldsymbol{\upphi}_{i,H}.
\end{eqnarray}

\begin{remark}
In (\ref{fH}) and (\ref{gdef}), the density function is the square sum of wavefunctions,
so we can expand (\ref{fH}) and (\ref{gdef}) into three parts, respectively.
Each term can be assembled efficiently.

However, this good property can not be used to exchange-correlation potential $V_{xc}$
because the structure of $V_{xc}$ is non-polynomial. For instance, the simplest form $V_{xc}(\rho)=-(\frac{3}{\pi}\rho)^{1/3}=-(\frac{3}{\pi}\sum_{i=1}^N(\phi_{i,H}+\theta_{2,i}\widetilde\phi_{i,h})^2)^{1/3}$,
then we can not expand it into different terms due to the fractional power $1/3$. This is also a reason to use the outer iteration to deal with  exchange-correlation potential.
Then the exchange-correlation potential $V_{xc}$ remains unchanged during the inner iteration.
\end{remark}

In Algorithm \ref{SCF_Strategy}, based on above discussion and preparation, we define the efficient
implementation strategy for Algorithm \ref{SCF_Iteration_Hartree}.
\begin{algorithm}[htbp]
\caption{Implementation strategy for Algorithm \ref{SCF_Iteration_Hartree}}\label{SCF_Strategy}
\begin{enumerate}
\item Preparation for the nonlinear iteration by assembling the following values for each linear boundary value problem:\\
Compute the matrices $A_{H,1}$, $A_{H,3}$, $A_{H,4}$ as in (\ref{Adef}), $A_{h4}$ as in (\ref{Ah4def}).\\
Compute the vectors $b_{H,1}$, $b_{H,3}$ and  $b_{H,4}$ as in (\ref{bdef}), ${\tt rhs}$ as in (\ref{rhsdef}).\\
Compute the scalars $\beta_1$, $\beta_3$ and $\beta_4$ as in (\ref{betadef}).\\
Compute the vector $g_{i,3}$ as in (\ref{gdef}).\\
Compute the matrix $G$ and vector $F_i$ as in (\ref{FG}). \\
\item Nonlinear iteration:
\begin{enumerate}
\item Compute the matrix $A_{H,2}$ as in (\ref{Adef}). Then
 $A_H=A_{H,1} + A_{H,2} + \theta_1A_{H,3}+ A_{H,4}.$
\item Compute the vector $b_{H,2}$ as in (\ref{bdef}). Then the vector $b_{Hh}=b_{H,1} + b_{H,2}+ \theta_1b_{H,3}+ b_{H,4}$.
\item Compute the scalar $\beta=\beta_1 + {\tt rhs}^T{\bf w_H} + \theta_1\beta_3+ \beta_4$.
\item Then solve the eigenvalue problem (\ref{Eigenvalue_H_h}) to get a new
eigenfunction $(\phi_{i,H},\theta_{2,i})$ and the corresponding eigenvalue $\lambda_{i,h}$.
\item Compute the vector $f_{H,i,1}$ as in (\ref{fH}). Then compute the vector $f_H$ according to (\ref{fH}).
\item Compute $g=\sum_{i=1}^N\big(\boldsymbol{\upphi}_{i,H}^TG \boldsymbol{\upphi}_{i,H}  + \theta_{2,i}F_i^T\boldsymbol{\upphi}_{i,H}+ g_{i,3}\big)$.
\item Then solve the boundary value problem (\ref{Hartree_H_h}) to get a new Hartree potential $w_H$ and $\theta_1$.
\item If the given accuracy is satisfied, stop the nonlinear iteration.
Otherwise, go to step (a) and continue the  nonlinear iteration.
\end{enumerate}
\item Output the eigenfunction $\phi_{i,h} = \phi_{i,H} + \theta_{2,i} \widetilde \phi_{i,h}$ and the eigenvalue $\lambda_{i,h}$.
\end{enumerate}
\end{algorithm}

\begin{remark}
We need to solve a Poisson equation to derive the Hartree potential. However, the Hartree potential does not decay exponentially as
the wavefunction. Thus, we can use the multipole expansion to derive a proper Dirichlet boundary condition. The detailed description
can be found in \cite{Baogang}, etc.
To further improve the solving efficiency, in our algorithm, the Dirichlet boundary condition will renew after each outer iteration
and remains unchanged in the inner iteration.
\end{remark}

In the last of this subsection, we want to emphasize that Algorithm \ref{Parallel_Aug_Subspace_Method} is naturally suitable for
eigenpairwise parallel computing due to its special structure.
As we can see, Algorithm \ref{Parallel_Aug_Subspace_Method} treats the different orbits independently.
This scheme can avoid doing inner products for orthogonalization in the high dimensional space $V_{h_{k+1}}$, which is always time-consuming and becomes
the bottle neck for parallel computing.
Thus parallel computing will benefit from such a strategy through solving different orbits in different processors.
During the iteration of Algorithm  \ref{Parallel_Aug_Subspace_Method}, we only need to transfer the data when it comes to compute
the density function, which just needs the wavefunction derived in each processor.
So the time spent on data transmission accounts for only a small part of the total computational time.
This is why Algorithm \ref{Parallel_Aug_Subspace_Method} is naturally suitable for parallel computing.

\section{Numerical results}
In this section, we provide several numerical examples to validate the
efficiency and scalability of the proposed numerical method
in this paper.  The numerical examples are carried out on LSSC-IV in
the State Key Laboratory of Scientific and Engineering Computing,
Chinese Academy of Sciences. Each computing node has two 18-core Intel
Xeon Gold 6140 processors at 2.3 GHz and 192 GB memory.  For the
involved small-scale eigenvalue problems, we adopt the implicitly
restarted Lanczos method provided in the package ARPACK \cite{arpack}.

The numerical experiments in this section are implemented for simulating the following five models.
The first model is the lithium hydroxide molecule and the complete Kohn-Sham equation is given by
\begin{equation}\label{13}
\left\{
\begin{array}{rcl}
\big(-\frac{1}{2}\Delta -\frac{3}{|x-r_1|}-\frac{1}{|x-r_2|}+
\int_{\Omega}\frac{\rho(y)}{|x-y|}dy+V_{xc}\big)\phi_i&=& \lambda_i\phi_i, \ \ \text{in } \Omega,\  \  i=1,2,\\
\phi_i&=& 0,\ \ \ \ \text{on } {\partial\Omega}, \ i=1,2,
\end{array}
\right.
\end{equation}
where $\Omega=(-6,6)^3$. In this equation, $\int_{\Omega}|\phi_i|^2dx=1 \ (i=1,2)$ and electron density
$\rho = 2(|\phi_1|^2+|\phi_2|^2)$, $r_1=(-1.0075,0,0)$, $r_2=(2.0075,0,0)$ denote the position of lithium atom and hydrogen atom,
and exchange-correlation potential is adopted as $V_{xc}(\rho)=-\frac{3}{2}\alpha (\frac{3}{\pi}\rho)^{1/3}$ with $\alpha = 0.77298$.
Since we don't consider spin polarisation, $2$ eigenpairs need to be calculated for this model.

The second model is Methane (CH4) and the computing domain is set to be $\Omega=(-6,6)^3$.  The associated parameters in the hamiltonian are the same as the first model.
For a full potential calculation, there are total $10$ electrons. Since we don't consider spin polarisation, $5$ eigenpairs need to be calculated.

The third model is Acetylene molecule (C2H2). We set the computing domain to be $\Omega=(-8,8)^3$.
For a full potential calculation, there are total fourteen electrons. Since we don't consider spin polarisation, seven eigenpairs need to be calculated.
In this model, the local density approximation (LDA) is defined as follows.
The exchange energy density is chosen as:
\begin{eqnarray*}
\epsilon_x(\rho)=-\frac{3}{4}\left(\frac{3}{\pi}\right)^{1/3}\rho(r)^{1/2}.
\end{eqnarray*}
The correlation energy density is chosen as \cite{perdew}:
\begin{equation*}
\epsilon_c(r_s)=\left\{
\begin{array}{rcl}
\ 0.0311  \text{ln} r_s-0.048+0.0020r_s\text{ln} r_s -0.0116 r_s,\ \ \ \ \text{if}\ \  r_s<1,\\
-0.1423/(1+1.0529\sqrt{r_s}+0.3334r_s),\ \ \ \text{if}\ \  r_s\geq 1
\end{array}
\right.
\end{equation*}
with $r_s=(\frac{3}{4\pi\rho})^{1/3}$.
The exchange-correlation potential is then chosen as
\begin{eqnarray*}
V_{xc}(\rho)=\epsilon_{xc}(\rho)+\rho\frac{d\epsilon_{xc}(\rho)}{d\rho}.
\end{eqnarray*}
In this model, there are total $14$ electrons and we need to compute $7$ eigenpairs for the full potential calculation.

The fourth model is Benzene molecule (C6H6) and the computing domain is $\Omega=(-10,10)^3$.
The corresponding parameters are the same as the third model with the LDA exchange-correlation potential.
There are total $42$ electrons. Since we don't consider spin polarisation,
$21$ eigenpairs need to be calculated for the full potential calculation.

The last model is the body-centered cubic (bcc) structure of Sodium crystal of $27$ cubes with $91$ atoms.
Here, we use the pseudopotential which was introduced in \cite{Nishioka} in our numerical experiment
and hence $91$ valence electrons and $46$ eigenpairs are simulated.
We use the LDA exchange-correlation potential as that in the third model.

Numerical experiments for these five molecules are implemented to
demonstrate following aspects, i.e., the feasibility of spearating
nonlinear Hartree potential and exchange-correlation potential in the
algorithm, the scalability of the algorithm, as well as the ability of
the algorithm on preserving the orthogonality etc.
In all numerical examples here,  the coarsest mesh $\mathcal T_H$ for  Algorithm  \ref{Parallel_Aug_Subspace_Method} is generated by the uniform refinement and
includes $24576$ mesh elements.

\subsection{On the feasibility}\label{Nonlinearity}
As mentioned in Introduction, there is significant difference on the
contribution for the total ground state energy from Hartree part and
exchange-correlation part. Hence, their performance in the convergence
towards the ground state should be also different. This can be
confirmed by following test.

We consider a lithium hydroxide molecule. To check the convergence
behavior with only exchange-correlation potential, governing equation
is given below,
\begin{equation}\label{11}
\left\{
\begin{array}{rcl}
\big(-\frac{1}{2}\Delta -\frac{3}{|x-r_1|}-\frac{1}{|x-r_2|}+V_{xc}(\rho)\big)\phi_i&=& \lambda_i\phi_i, \ \ \text{in } \Omega,\  \  i=1,2,\\
\phi_i&=& 0,\ \ \ \ \text{on } {\partial\Omega}, \ i=1,2,
\end{array}
\right.
\end{equation}
where $\Omega$, $\rho$, $r_1$, $r_2$ and $V_{xc}(\rho)$ are the same as in (\ref{13}).


We use the standard AFEM to solve (\ref{11}) with the refinement parameter $\theta=0.4$. Table \ref{ex1-table1} presents
the number of SCF iteration times in each level of the adaptive finite element spaces.

For the comparison, we next consider the same molecule, but with only
the Hartree potential in the hamiltonian,
\begin{equation}\label{12}
\left\{
\begin{array}{rcl}
\big(-\frac{1}{2}\Delta -\frac{3}{|x-r_1|}-\frac{1}{|x-r_2|}+
\int_{\Omega}\frac{\rho(y)}{|x-y|}dy\big)\phi_i&=& \lambda_i\phi_i, \ \ \text{in } \Omega,\  \  i=1,2,\\
\phi_i&=& 0,\ \ \ \ \text{on } {\partial\Omega}, \ i=1,2.
\end{array}
\right.
\end{equation}
Table \ref{ex1-table2} presents the corresponding number of SCF iteration times in each level of adaptive finite element spaces by the standard AFEM.
\begin{table}[htbp]
\begin{center}
\begin{tabular}{|c|c|c|c|c|c|c|c|c|c|c|c|c|c|c|c|c|c|}\hline
Mesh Level &1& 2 &3 &4&5&6& 7 &8 &9&10&11& 12 &13 &14&15  \\ \hline
    SCF iteration&07& 07 & 07&07&06 &07& 07 & 07 &07&06&07& 07 &07 & 06&07 \\ \hline
\end{tabular}
\end{center}
\caption{The SCF iteration numbers for equation (\ref{11}).}
\label{ex1-table1}
\end{table}

\begin{table}[htbp]
\begin{center}
\begin{tabular}{|c|c|c|c|c|c|c|c|c|c|c|c|c|c|c|c|c|c|}\hline
Mesh Level &1& 2 &3 &4&5&6& 7 &8 &9&10&11& 12 &13 &14&15  \\ \hline
    SCF iteration&23& 22 & 23&22&22 &22& 23 & 22 &23&22&22&22 &22 &22&22 \\ \hline
\end{tabular}
\end{center}
\caption{The SCF iteration numbers for  equation (\ref{12}).}
\label{ex1-table2}
\end{table}

From Tables \ref{ex1-table1} and \ref{ex1-table2}, it can be observed
clearly that Hartree potential brings more iterations towards the
ground state, compared with exchange-correlation potential. A similar
observation is always available for other molecules.



For further understanding the relation of the outer and inner iterations, we solve
the lithium hydroxide molecule which is defined by (\ref{13}).
Table \ref{ex1-table3} presents the numbers of the inner iterations
and outer iterations of Algorithm \ref{Parallel_Aug_Subspace_Method}
in the final five adaptive finite element spaces.  Because the approximate solution obtained
from the last outer iteration is used as the initial value, we can see
that with approaching the ground state of the molecule, the number of
the inner iteration becomes smaller and smaller. Together with the
acceleration strategy proposed in Subsection \ref{KS_CH4}, the efficiency
of proposed algorithm will be improved significantly, compared with
the multilevel correction method proposed in \cite{HuXieXu}. The
comparison of the efficiency will be given in following subsections.

\begin{table}[htbp]
\begin{center}
\begin{tabular}{|c|c|c|c|c|c|c|c|c|}\hline
Outer iteration &1& 2 &3 &4&5&6&7&8  \\ \hline
Inner iteration &22& 17 & 12 &8&5&3&2&0 \\ \hline
\end{tabular}

\begin{tabular}{|c|c|c|c|c|c|c|c|c|}\hline
Outer iteration &1& 2 &3 &4&5&6&7&8  \\ \hline
Inner iteration &20& 15 & 11 &8&5&2&2&0 \\ \hline
\end{tabular}

\begin{tabular}{|c|c|c|c|c|c|c|c|c|}\hline
Outer iteration &1& 2 &3 &4&5&6&7&8  \\ \hline
Inner iteration &21& 15 & 12 &8&6&3&2&1 \\ \hline
\end{tabular}

\begin{tabular}{|c|c|c|c|c|c|c|c|c|}\hline
Outer iteration &1& 2 &3 &4&5&6&7&8  \\ \hline
Inner iteration &22& 16 & 10 &6&3&2&1&0 \\ \hline
\end{tabular}

\begin{tabular}{|c|c|c|c|c|c|c|c|c|}\hline
Outer iteration &1& 2 &3 &4&5&6&7&8  \\ \hline
Inner iteration &22& 16 & 11 &8&4&2&1&0 \\ \hline
\end{tabular}

\end{center}
\caption{The outer and inner iteration numbers in the final five adaptive finite element spaces for Hydrogen-Lithium.}
\label{ex1-table3}
\end{table}

%

\begin{table}[htbp]
\begin{center}
\begin{tabular}{|c|c|c|c|c|c|c|c|c|}\hline
Outer iteration &1& 2 &3 &4&5&6&7&8  \\ \hline
Inner iteration &32& 24 & 15 &12&9&6&3&2 \\ \hline
\end{tabular}

\begin{tabular}{|c|c|c|c|c|c|c|c|c|}\hline
Outer iteration &1& 2 &3 &4&5&6&7&8  \\ \hline
Inner iteration &33& 25 & 18 &12&8&5&3&3 \\ \hline
\end{tabular}

\begin{tabular}{|c|c|c|c|c|c|c|c|c|}\hline
Outer iteration &1& 2 &3 &4&5&6&7&8  \\ \hline
Inner iteration &32& 25 & 17 &11&8&5&3&1 \\ \hline
\end{tabular}

\begin{tabular}{|c|c|c|c|c|c|c|c|c|}\hline
Outer iteration &1& 2 &3 &4&5&6&7&8  \\ \hline
Inner iteration &32& 25 & 17 &12&8&6&3&2 \\ \hline
\end{tabular}

\begin{tabular}{|c|c|c|c|c|c|c|c|c|}\hline
Outer iteration &1& 2 &3 &4&5&6&7&8  \\ \hline
Inner iteration &33& 26 & 19 &14&9&6&3&2 \\ \hline
\end{tabular}
\end{center}
\caption{The outer and inner iteration numbers in the final five adaptive finite element spaces for Methane.}
\label{ex2-table1}
\end{table}

\begin{table}[htbp]
\begin{center}
\begin{tabular}{|c|c|c|c|c|c|c|c|}\hline
Outer iteration &1& 2 &3 &4&5&6&7 \\ \hline
Inner iteration &50& 37 & 26 &16&08&3&1 \\ \hline
\end{tabular}

\begin{tabular}{|c|c|c|c|c|c|c|c|}\hline
Outer iteration &1& 2 &3 &4&5&6&7  \\ \hline
Inner iteration &53& 42 & 30 &19&11&5&2 \\ \hline
\end{tabular}

\begin{tabular}{|c|c|c|c|c|c|c|c|}\hline
Outer iteration &1& 2 &3 &4&5&6&7  \\ \hline
Inner iteration &50& 38 & 26 &15&06&2&0 \\ \hline
\end{tabular}

\begin{tabular}{|c|c|c|c|c|c|c|c|}\hline
Outer iteration &1& 2 &3 &4&5&6&7  \\ \hline
Inner iteration &52& 40 & 29 &19&10&3&1 \\ \hline
\end{tabular}

\begin{tabular}{|c|c|c|c|c|c|c|c|}\hline
Outer iteration &1& 2 &3 &4&5&6&7  \\ \hline
Inner iteration &50& 34 & 20 &11&05&3&0 \\ \hline
\end{tabular}
\end{center}
\caption{The outer and inner iteration numbers in the final five adaptive finite element spaces for Acetylene.}
\label{ex3-table2}
\end{table}

\begin{table}[htbp]
\begin{center}
\begin{tabular}{|c|c|c|c|c|c|c|c|c|c|c|c|}\hline
Outer iteration &1& 2 &3 &4&5&6&7&8&9&10&11  \\ \hline
Inner iteration &177& 135 &108 & 72 & 55 & 50 & 42 & 34 & 22 & 10 & 5  \\ \hline
\end{tabular}

\begin{tabular}{|c|c|c|c|c|c|c|c|c|c|c|c|}\hline
Outer iteration &1& 2 &3 &4&5&6&7&8&9&10&11  \\ \hline
Inner iteration &177& 134 & 095 & 67 & 50 & 38 & 26 & 16 & 08 & 3 & 0 \\ \hline
\end{tabular}

\begin{tabular}{|c|c|c|c|c|c|c|c|c|c|c|c|}\hline
Outer iteration &1& 2 &3 &4&5&6&7&8&9&10&11  \\ \hline
Inner iteration &165& 138 & 090 &64&45&34&25&17&10&05&0 \\ \hline
\end{tabular}

\begin{tabular}{|c|c|c|c|c|c|c|c|c|c|c|c|}\hline
Outer iteration &1& 2 &3 &4&5&6&7&8&9&10&11  \\ \hline
Inner iteration &175& 135 &103 &74&53&40&30&22&14&8&4 \\ \hline
\end{tabular}

\begin{tabular}{|c|c|c|c|c|c|c|c|c|c|c|c|}\hline
Outer iteration &1& 2 &3 &4&5&6&7&8&9&10&11  \\ \hline
Inner iteration &170& 141 &108 &81&61&47&37&26&18&9&4 \\ \hline
\end{tabular}
\end{center}
\caption{The outer and inner iteration numbers  in the final five adaptive finite element spaces for Benzene.}
\label{ex4-table1}
\end{table}

\begin{table}[htbp]
\begin{center}
\begin{tabular}{|c|c|c|c|c|c|c|c|c|c|c|c|c|}\hline
Outer iteration &1& 2 &3 &4&5&6&7&8&9&10&11&12  \\ \hline
Inner iteration &313& 248 &184 & 147 & 112 & 88 & 64 & 48 & 33 & 21 & 12 &4 \\ \hline
\end{tabular}

\begin{tabular}{|c|c|c|c|c|c|c|c|c|c|c|c|c|}\hline
Outer iteration &1& 2 &3 &4&5&6&7&8&9&10&11&12  \\ \hline
Inner iteration &312& 246 & 180 & 145 & 111 & 87 & 60 & 47 & 31 & 19 & 10 &3\\ \hline
\end{tabular}

\begin{tabular}{|c|c|c|c|c|c|c|c|c|c|c|c|c|}\hline
Outer iteration &1& 2 &3 &4&5&6&7&8&9&10&11&12  \\ \hline
Inner iteration &310& 243 & 173 &144&107&84&62&44&30&11&0&3 \\ \hline
\end{tabular}

\begin{tabular}{|c|c|c|c|c|c|c|c|c|c|c|c|c|}\hline
Outer iteration &1& 2 &3 &4&5&6&7&8&9&10&11&12  \\ \hline
Inner iteration &306& 238 &162 &136&102&80&61&40&26&13&5&0 \\ \hline
\end{tabular}

\begin{tabular}{|c|c|c|c|c|c|c|c|c|c|c|c|c|}\hline
Outer iteration &1& 2 &3 &4&5&6&7&8&9&10&11 &12 \\ \hline
Inner iteration &310& 240 &175 &146&105&87&63&42&29&18&11&4 \\ \hline
\end{tabular}

\end{center}
\caption{The outer and inner iteration numbers in the final five adaptive finite element spaces for Sodium crystal.}
\label{ex5-table1}
\end{table}

It should be pointed out that similar results are also available for
other molecules. Tables \ref{ex2-table1}-\ref{ex5-table1} present the
results of Algorithm \ref{Parallel_Aug_Subspace_Method} in the final
five adaptive finite element spaces, for Methane, Acetylene, Benzene
and Sodium crystal molecules.
\subsection{Parallel scalability testing of Algorithm \ref{Parallel_Aug_Subspace_Method}}\label{scalability}
The aim of this subsection is to check the parallel scalability of the eigenpairwise strategy in Algorithm \ref{Parallel_Aug_Subspace_Method}.
In the step (a) of Algorithm  \ref{SCF_Iteration_Hartree}, we solve the $N$ linear boundary value problem in parallel since they
are independent from each other. In the step (b) of Algorithm  \ref{SCF_Iteration_Hartree}, the right hand side term of (\ref{hartree})
requires the eigenfunctions from the $N$ linear boundary value problems.
In our numerical experiments, each term in the right hand side term of (\ref{hartree}) will be computed in parallel on different computing nodes.
Then the data will be sent to one node to form the right hand side term of (\ref{hartree}).
Finally, the equation (\ref{hartree}) for Hartree potential is solved on this node, and then be broadcasted to all the computing nodes to perform the next loop.
During the numerical experiments, the communication between different computing nodes is realized through MPI.

\begin{figure}[htbp]
\centering
\includegraphics[width=5.5cm]{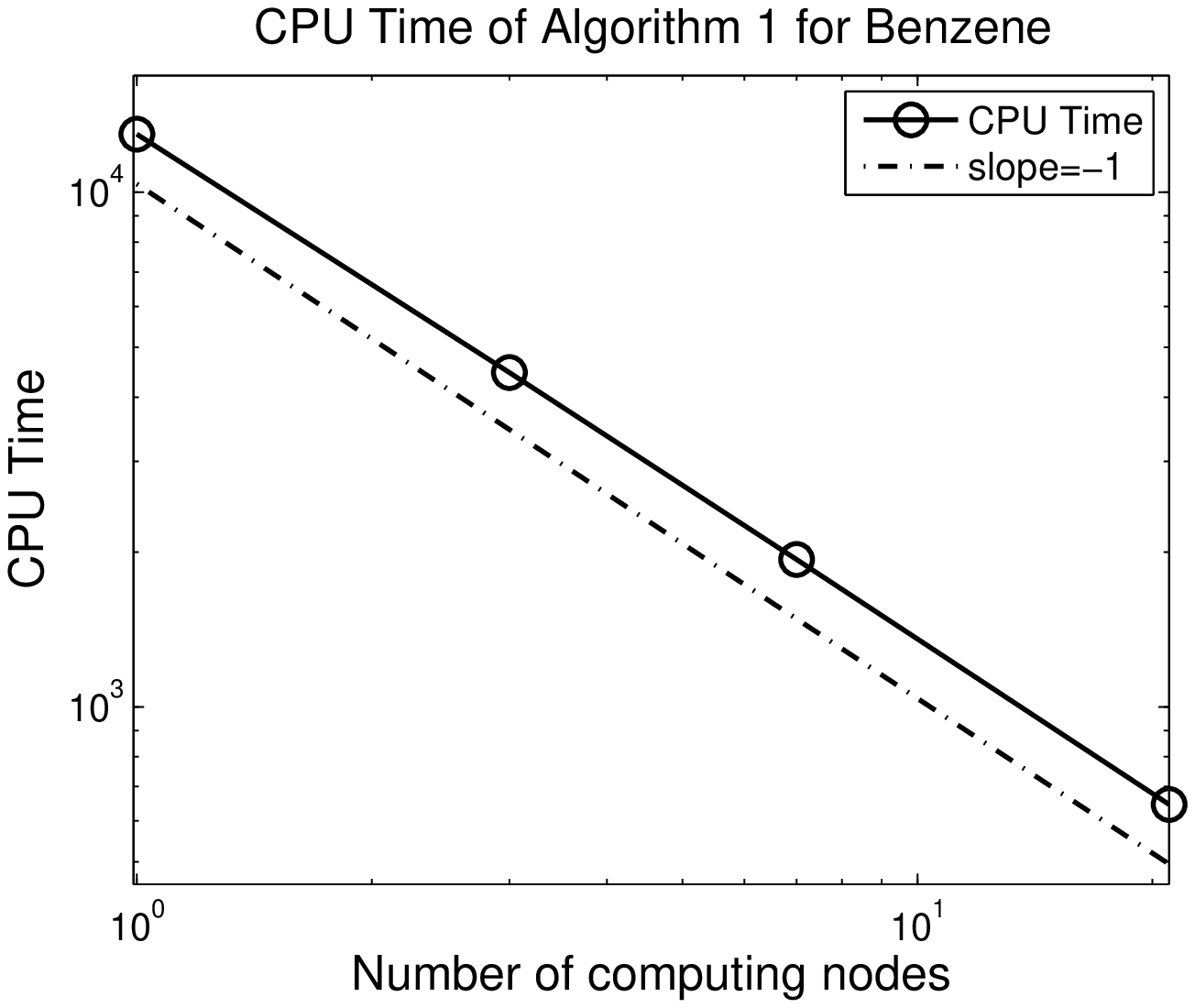}\ \ \ \ \ \ \ \ \
\includegraphics[width=5.5cm]{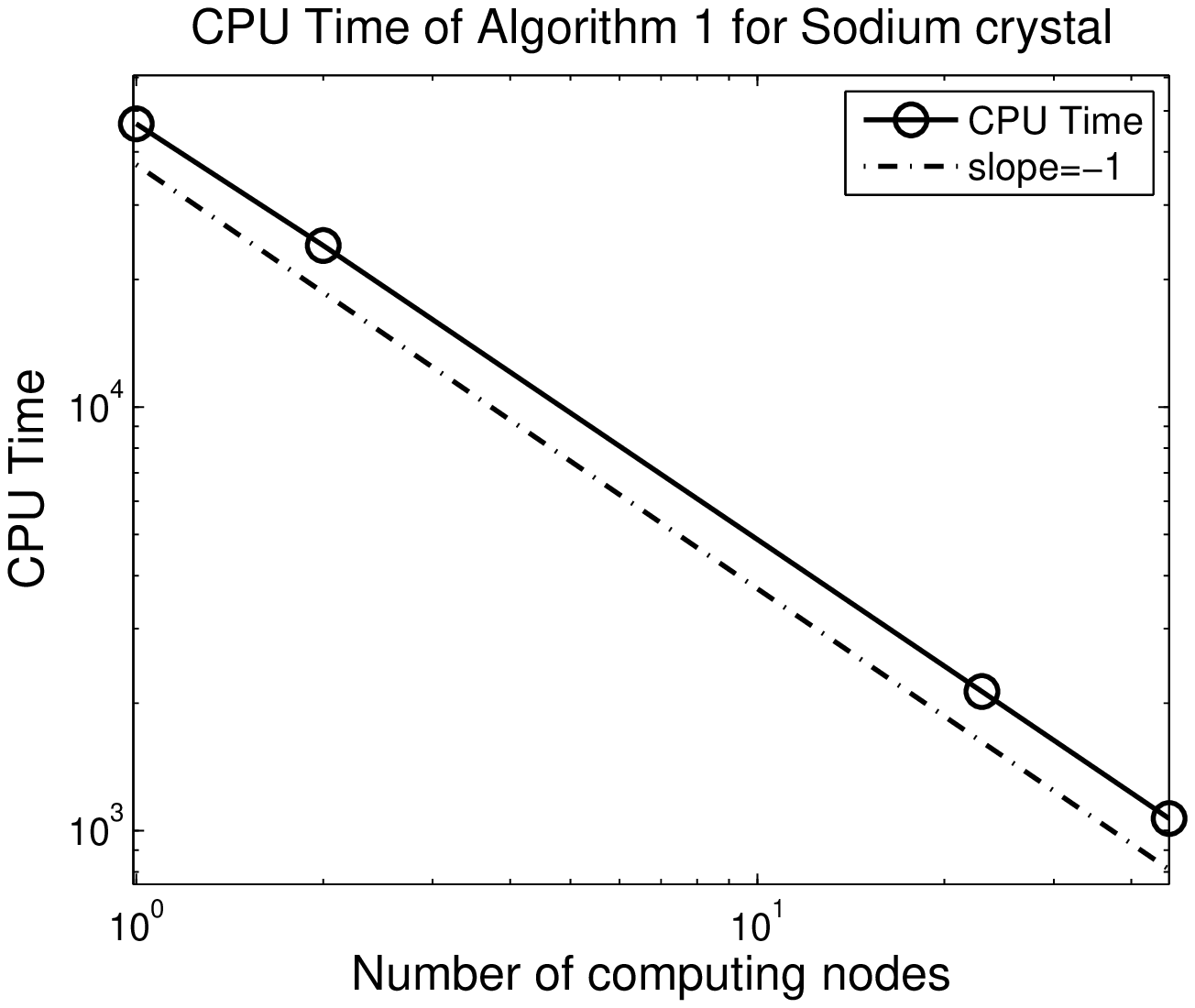}
\caption{Computational time of Algorithm \ref{Parallel_Aug_Subspace_Method} for Benzene (left) and Sodium crystal (right).}
\label{ben-1}
\end{figure}
\begin{figure}[htbp]
\centering
\includegraphics[width=5.5cm]{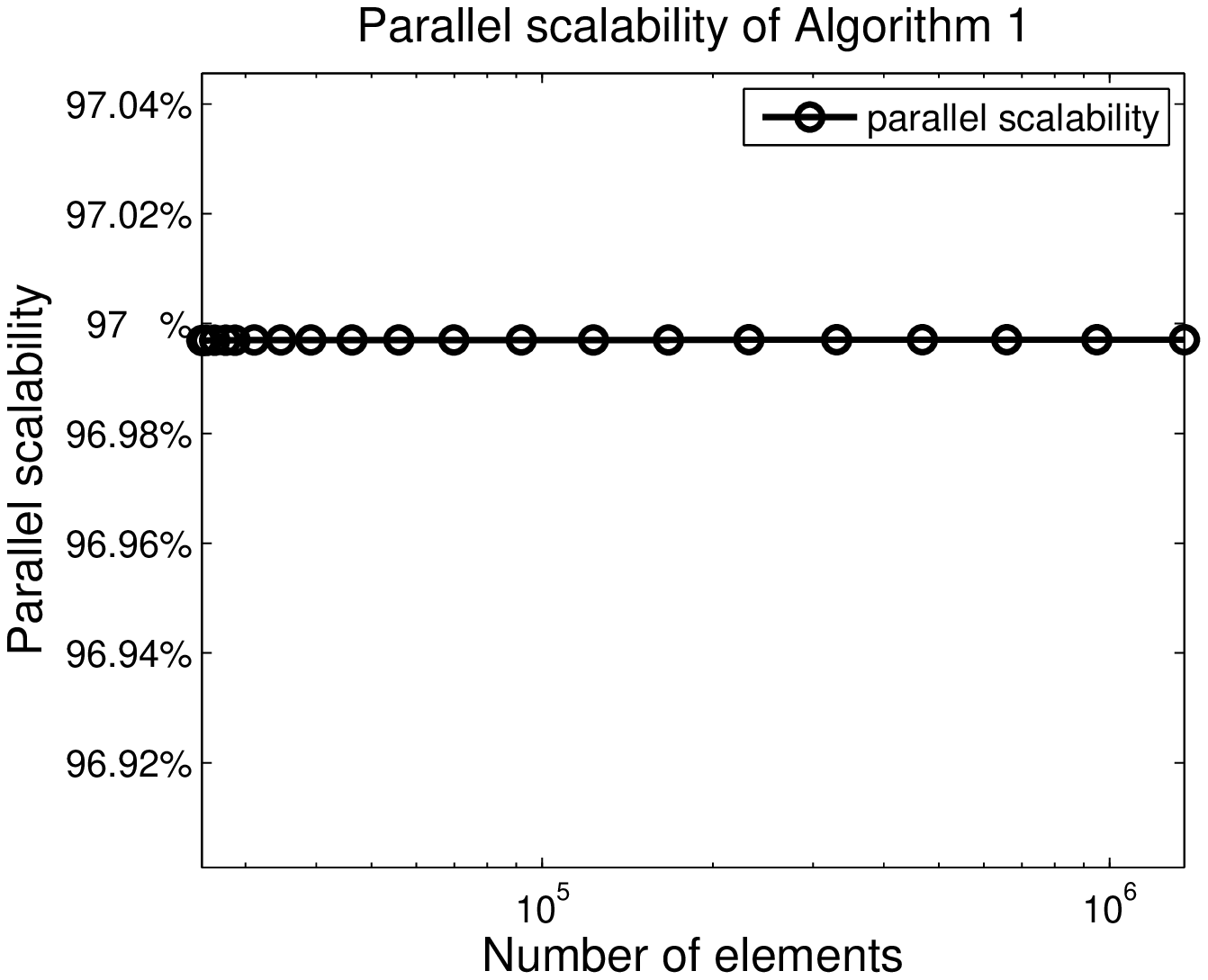}\ \ \ \ \ \ \ \ \
\includegraphics[width=5.5cm]{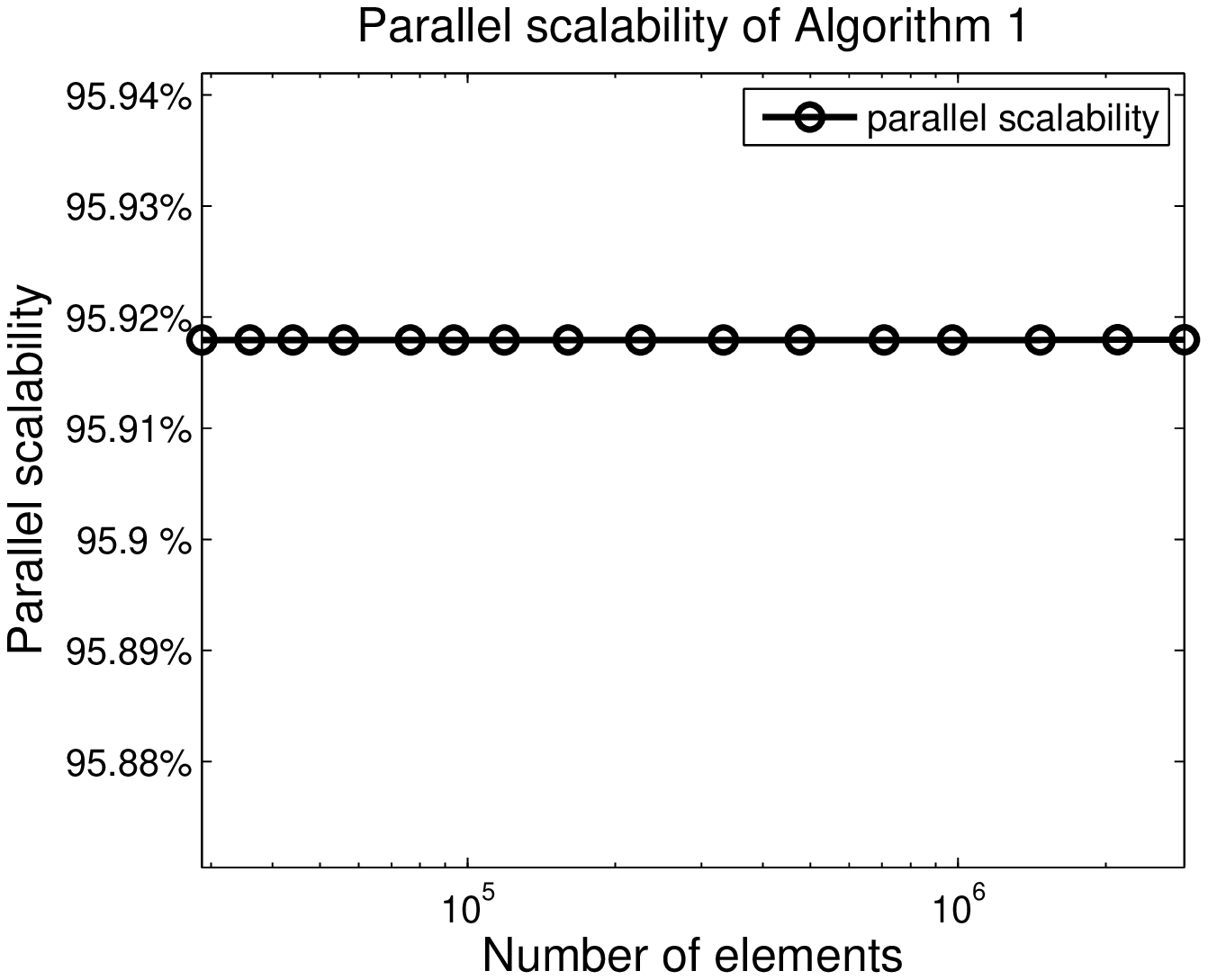}
\caption{Parallel efficiency of Algorithm \ref{Parallel_Aug_Subspace_Method} for Hydrogen-Lithium (left) and Methane (right).}
\label{para-1}
\end{figure}
\begin{figure}[htbp]
\centering
\includegraphics[width=5.5cm]{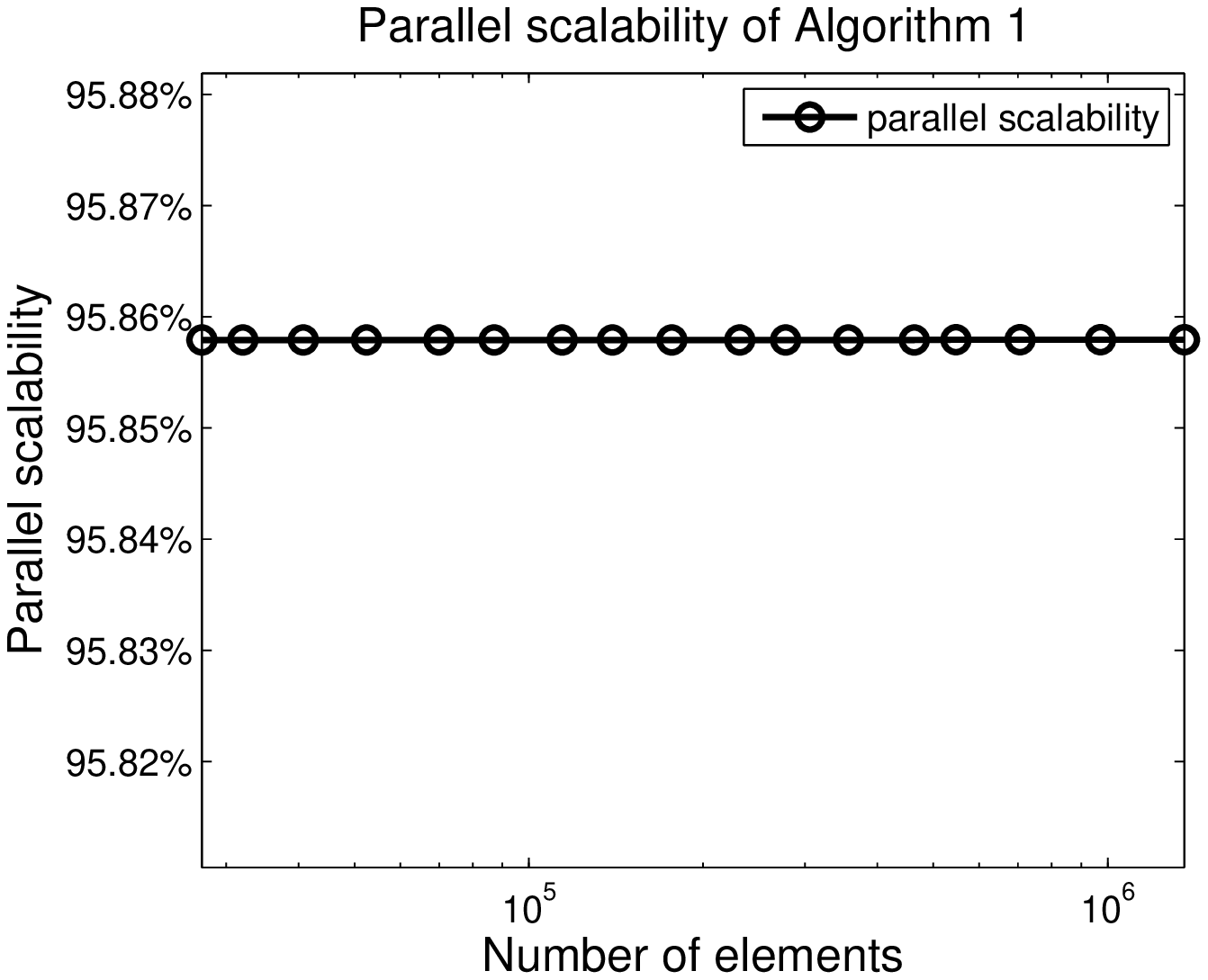}\ \ \ \
\includegraphics[width=5.5cm]{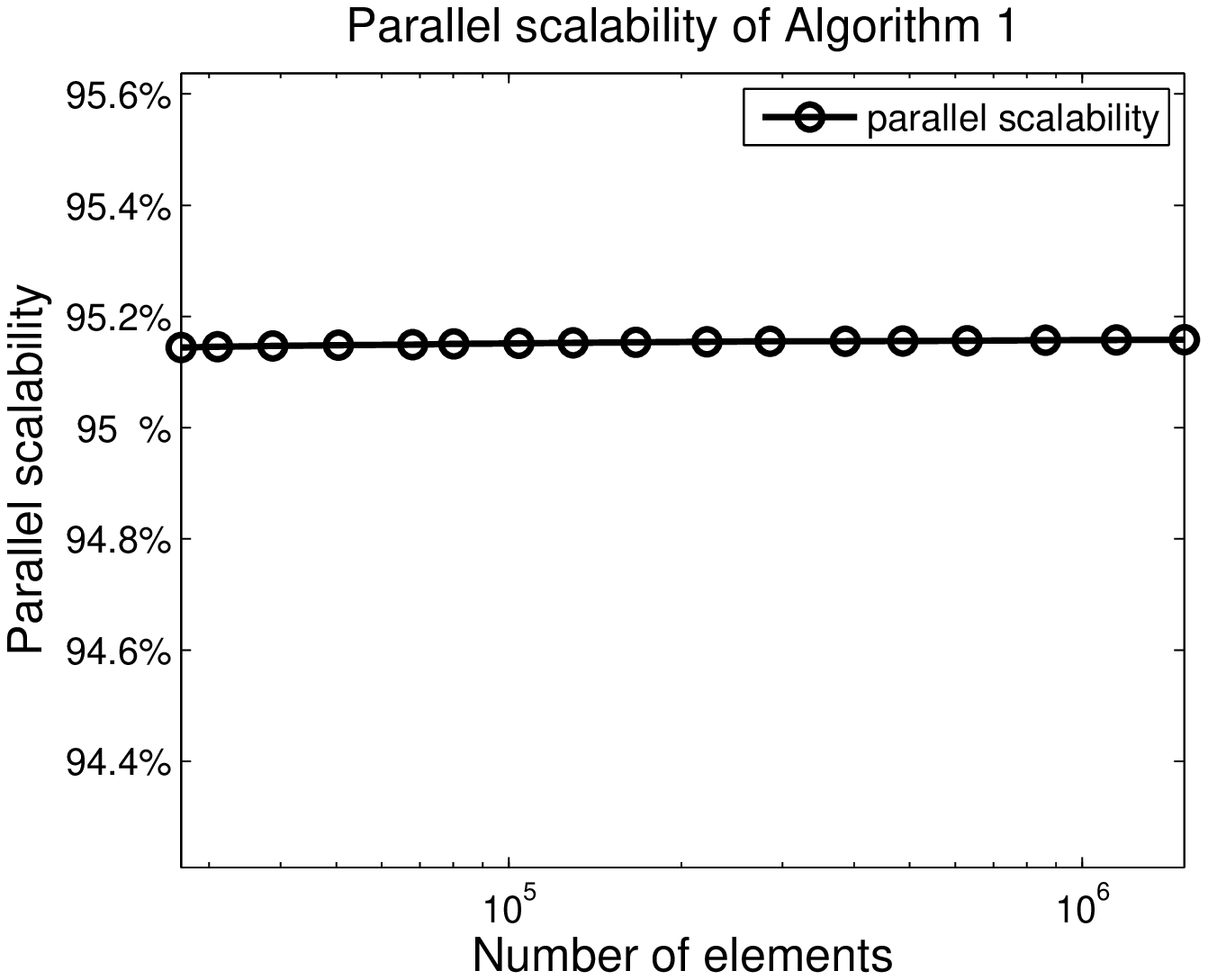}\ \ \ \
\includegraphics[width=5.5cm]{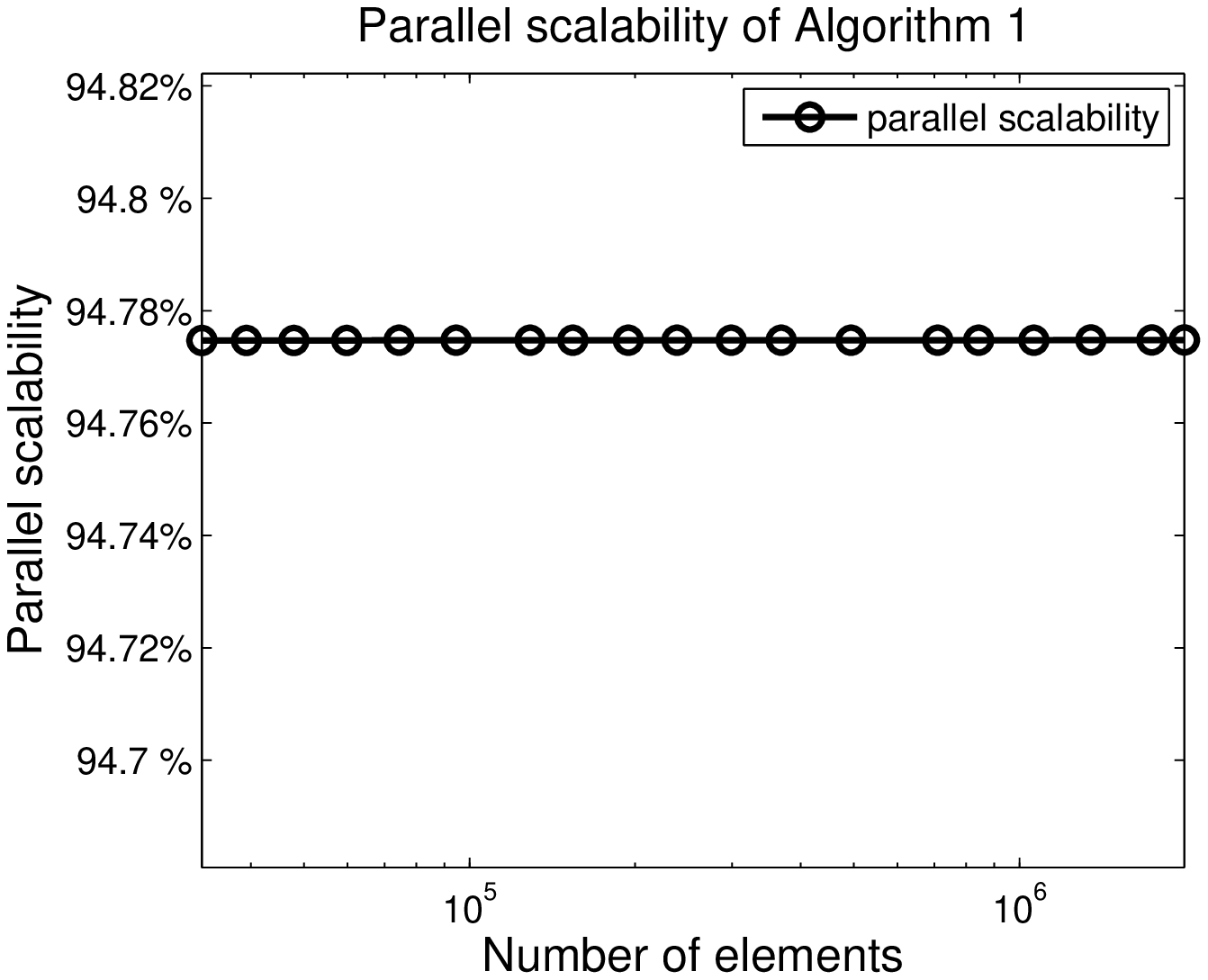}
\caption{Parallel efficiency of Algorithm \ref{Parallel_Aug_Subspace_Method} for Acetylene (left), Benzene (middle) and Sodium crystal (right).}
\label{para-2}
\end{figure}

In order to show the scalability of Algorithm \ref{Parallel_Aug_Subspace_Method}, we first present the numerical results
for Benzene and Sodium crystal.
We test the parallel scalability for Benzene and Sodium crystal when different numbers of computing nodes are used.
The desired eigenpairs are distributed to different computing nodes and
each computing node has the same number of eigepairs to be solved.
The number of computing nodes for Benzene and Sodium crystal are set to be $[1\ 3\ 7\ 21]$ and $[1\ 2\ 23\ 46]$, respectively.
The total computational time of Algorithm \ref{Parallel_Aug_Subspace_Method} by using different number of computing nodes
are presented in Figure \ref{ben-1}.
From Figure \ref{ben-1}, we can find that Algorithm \ref{Parallel_Aug_Subspace_Method} has a good scalability.

Besides, we also test the parallel efficiency of Algorithm \ref{Parallel_Aug_Subspace_Method}
when the mesh size changed while the number of computing nodes remain unchanged.
The parallel efficiency is calculated through the formula $(T_s/T_p)/N_p$,
where $T_s$ denote the serial computing time, $T_p$ denotes the parallel computing time and $N_p$ denotes the number of computing nodes.
The numbers of computing nodes are chosen as that of the required eigenpairs for the five models, i.e. $[2\ 5\ 7\ 21\ 46]$, respectively.
We demonstrate the trend of parallel efficiency according to the refinement of mesh and
the corresponding results are presented in Figures \ref{para-1}-\ref{para-2}.
From Figures \ref{para-1}-\ref{para-2}, we can also find that
Algorithm \ref{Parallel_Aug_Subspace_Method} has a good scalability in different adaptive spaces.

\subsection{Numerical performance of Algorithm \ref{Parallel_Aug_Subspace_Method}}\label{KS_CH4}
In this subsection, we show the numerical performances for the five
Kohn-Sham models.  We first test the orthogonality of the approximate
eigenfunctions derived from Algorithm
\ref{Parallel_Aug_Subspace_Method}.  Tables
\ref{ex1-table4}-\ref{ex5-table2} present the biggest values of inner
products of the eigenfunctions corresponding to the different
eigenvalues in each level of the adaptive finite element spaces for
the five Kohn-Sham models.  The results presented in Tables
\ref{ex1-table4}-\ref{ex5-table2} show that Algorithm
\ref{Parallel_Aug_Subspace_Method} can keep the orthogonality of the
eigenfunctions along with the refinement of mesh.

Next, we investigate the error estimates of the approximate solutions
derived by Algorithm \ref{Parallel_Aug_Subspace_Method} for the five models.

\begin{figure}[htbp]
\centering
\includegraphics[width=5.0cm]{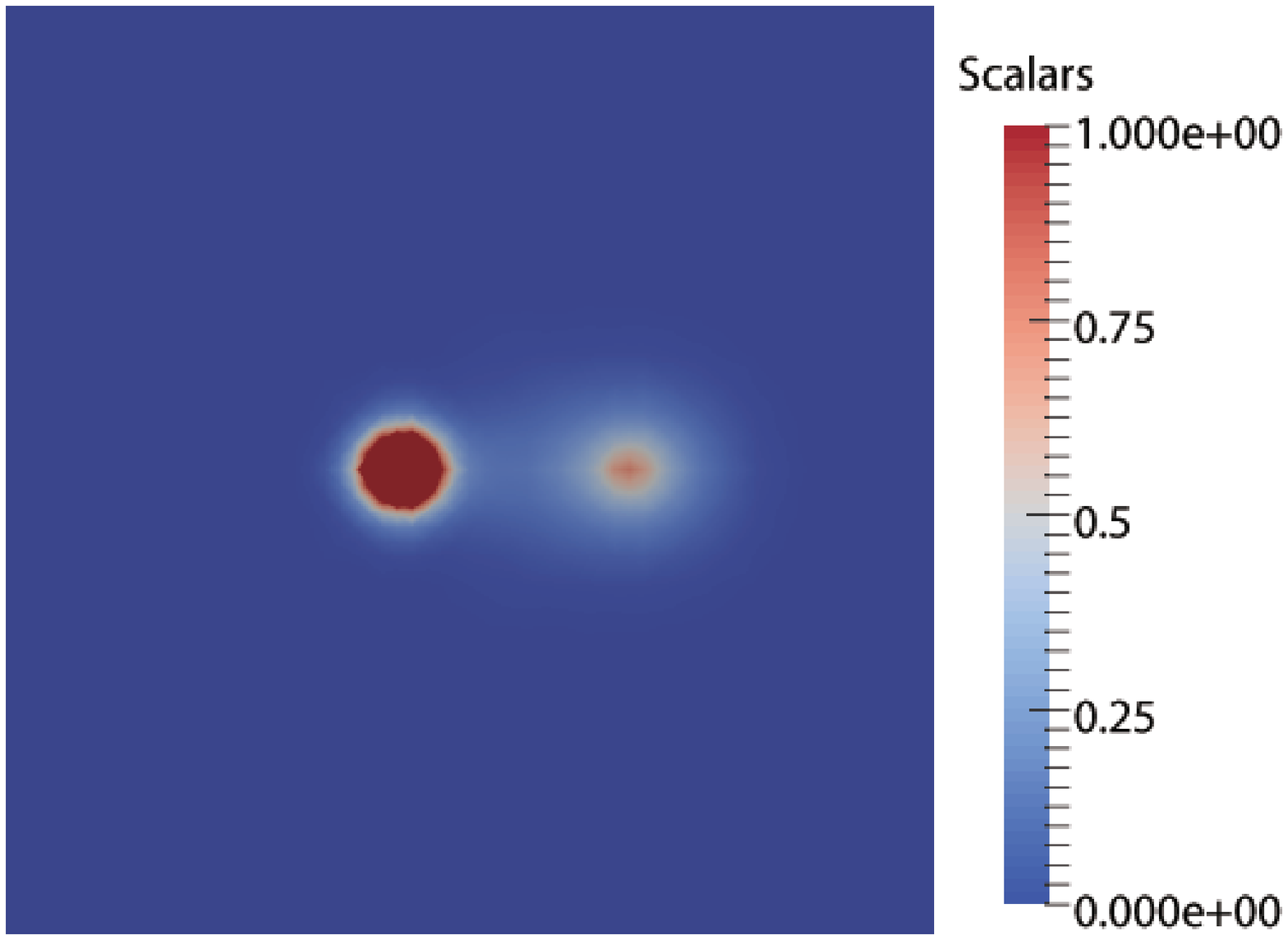}\ \ \ \ \ \ \ \ \ \ \ \ \
\includegraphics[width=3.5cm]{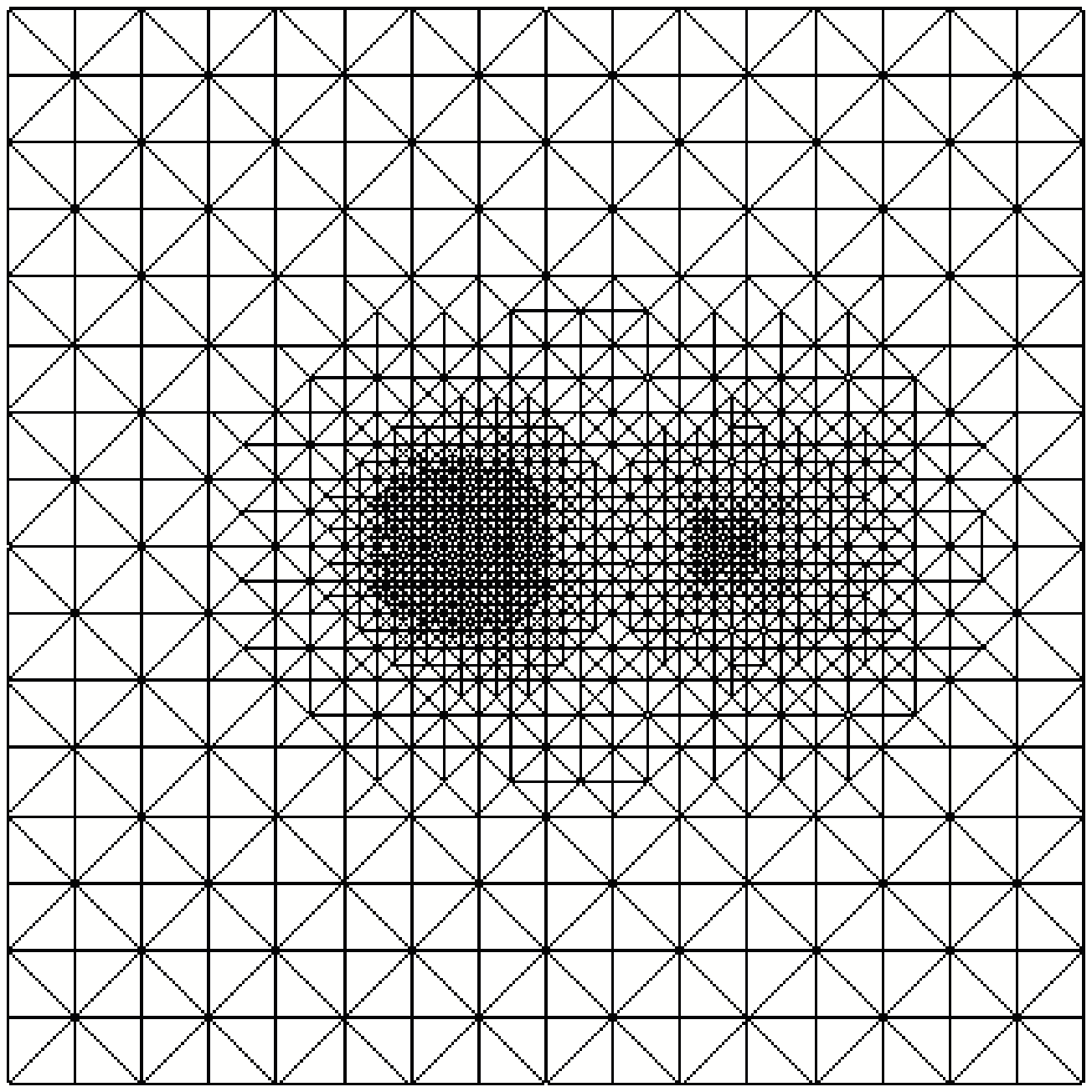}
\caption{Contour plot of the electron density (left) and the mesh after 10 adaptive refinements (righth) of Algorithm \ref{Parallel_Aug_Subspace_Method} for Hydrogen-Lithium.}
\label{ex1-f1}
\end{figure}

\begin{figure}[htbp]
\centering
\includegraphics[height=5.5cm, width=5.5cm]{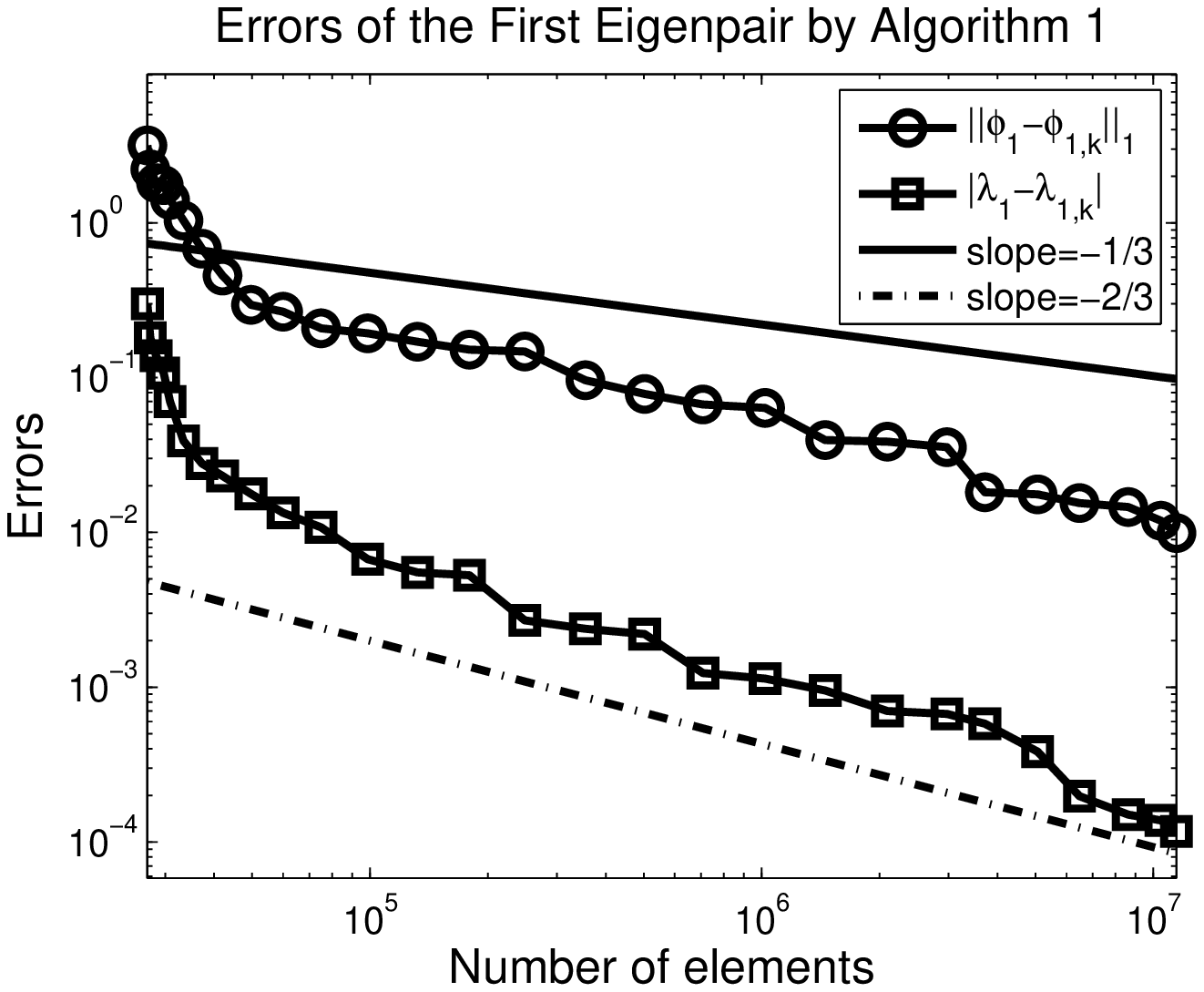}
\includegraphics[height=5.5cm, width=5.5cm]{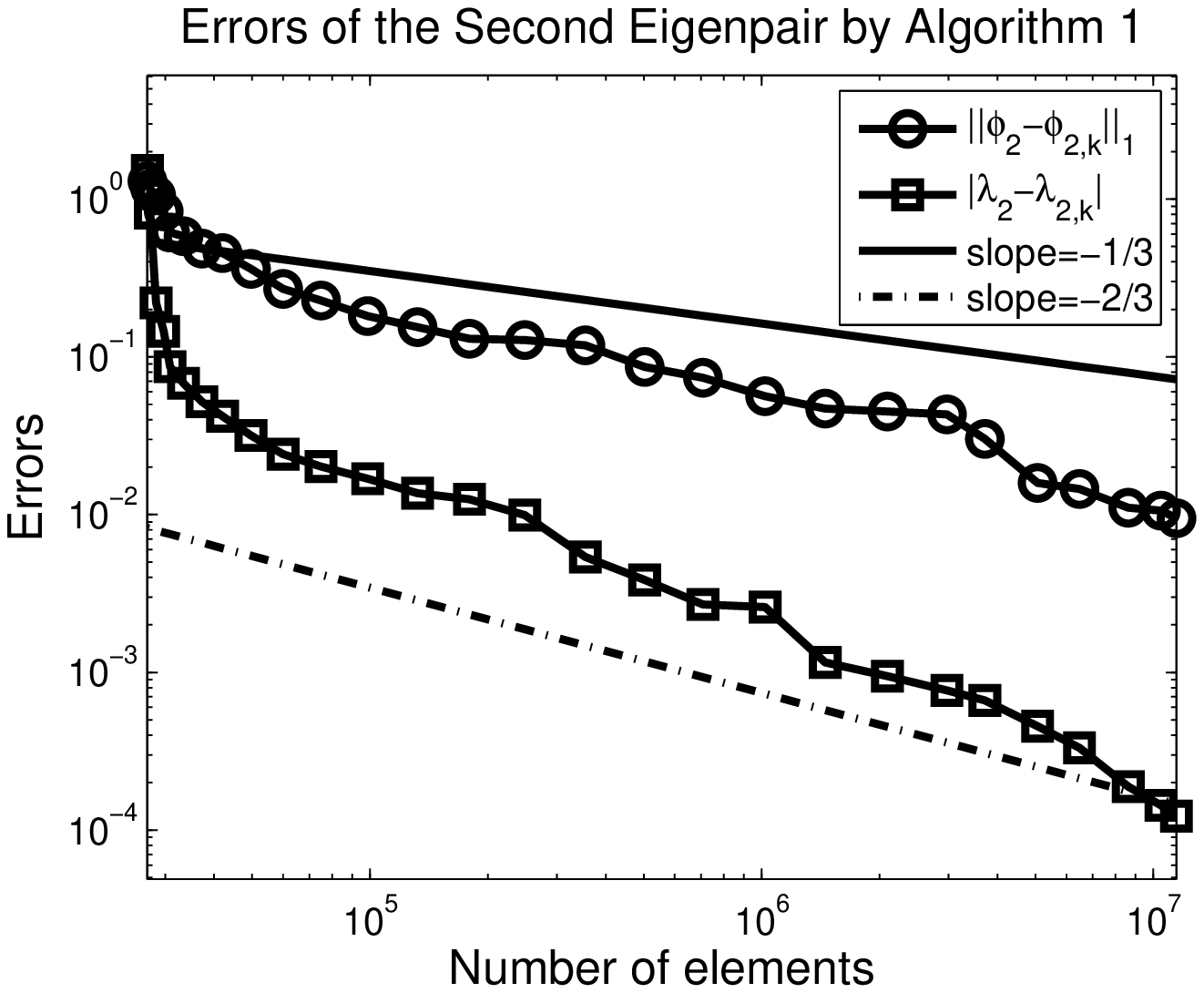}
\includegraphics[height=5.5cm, width=5.5cm]{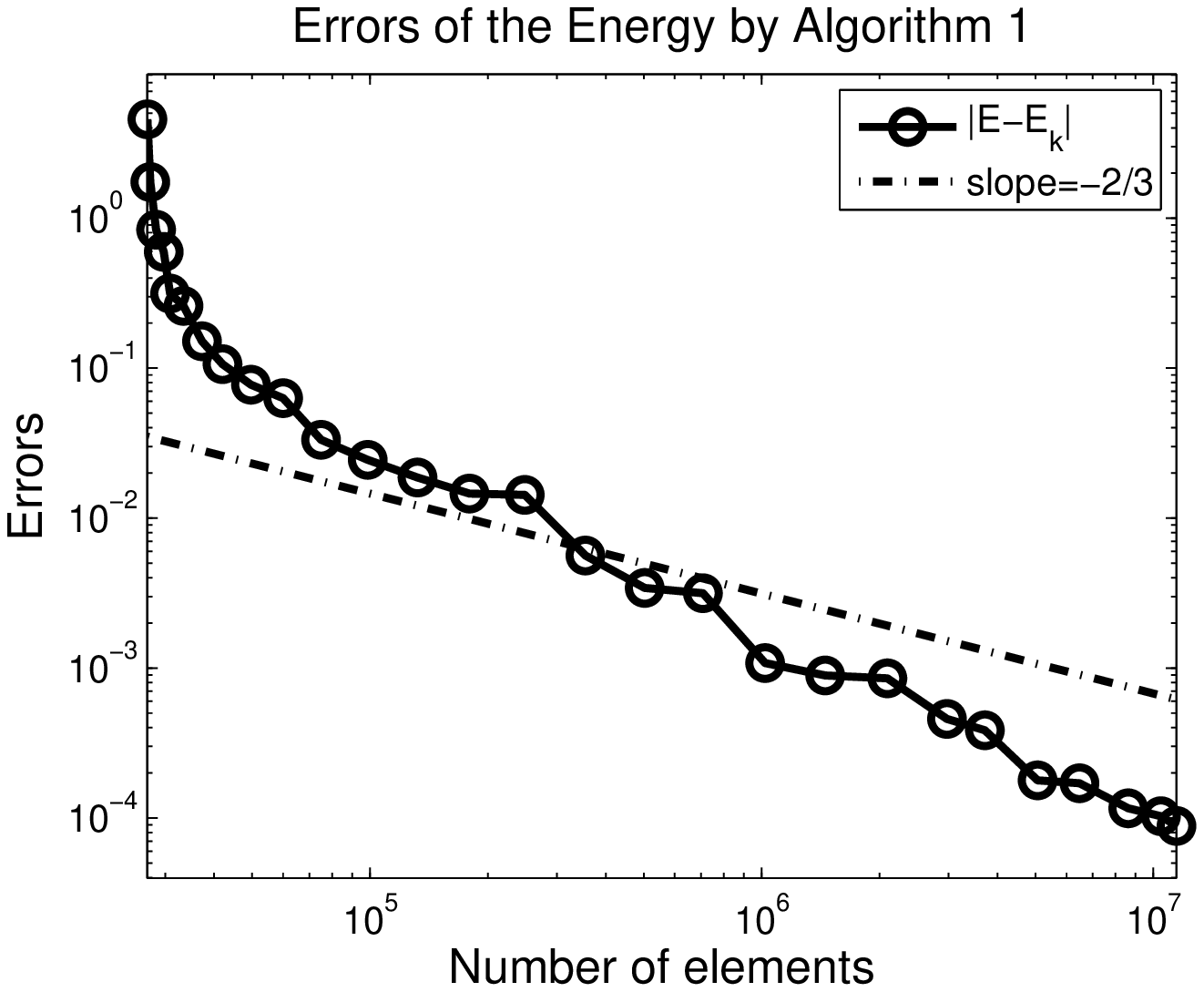}
\caption{Errors of eigenpair approximations (left and middle)
and energy approximations (right) of Algorithm \ref{Parallel_Aug_Subspace_Method} for Hydrogen-Lithium.}
\label{ex1-f2}
\end{figure}

\begin{figure}[htbp]
\centering
\includegraphics[width=3.7cm]{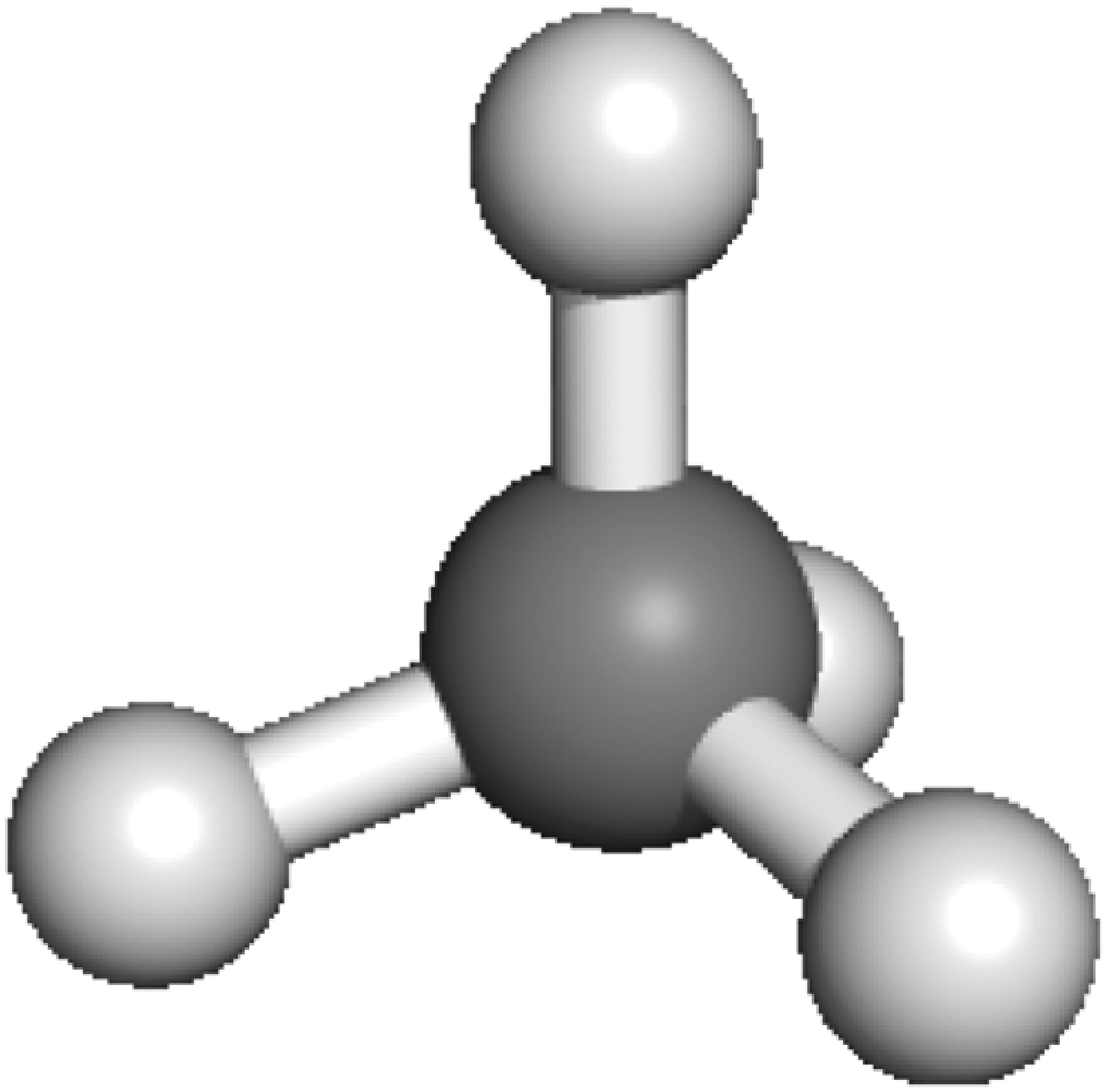}
 \includegraphics[width=3.6cm]{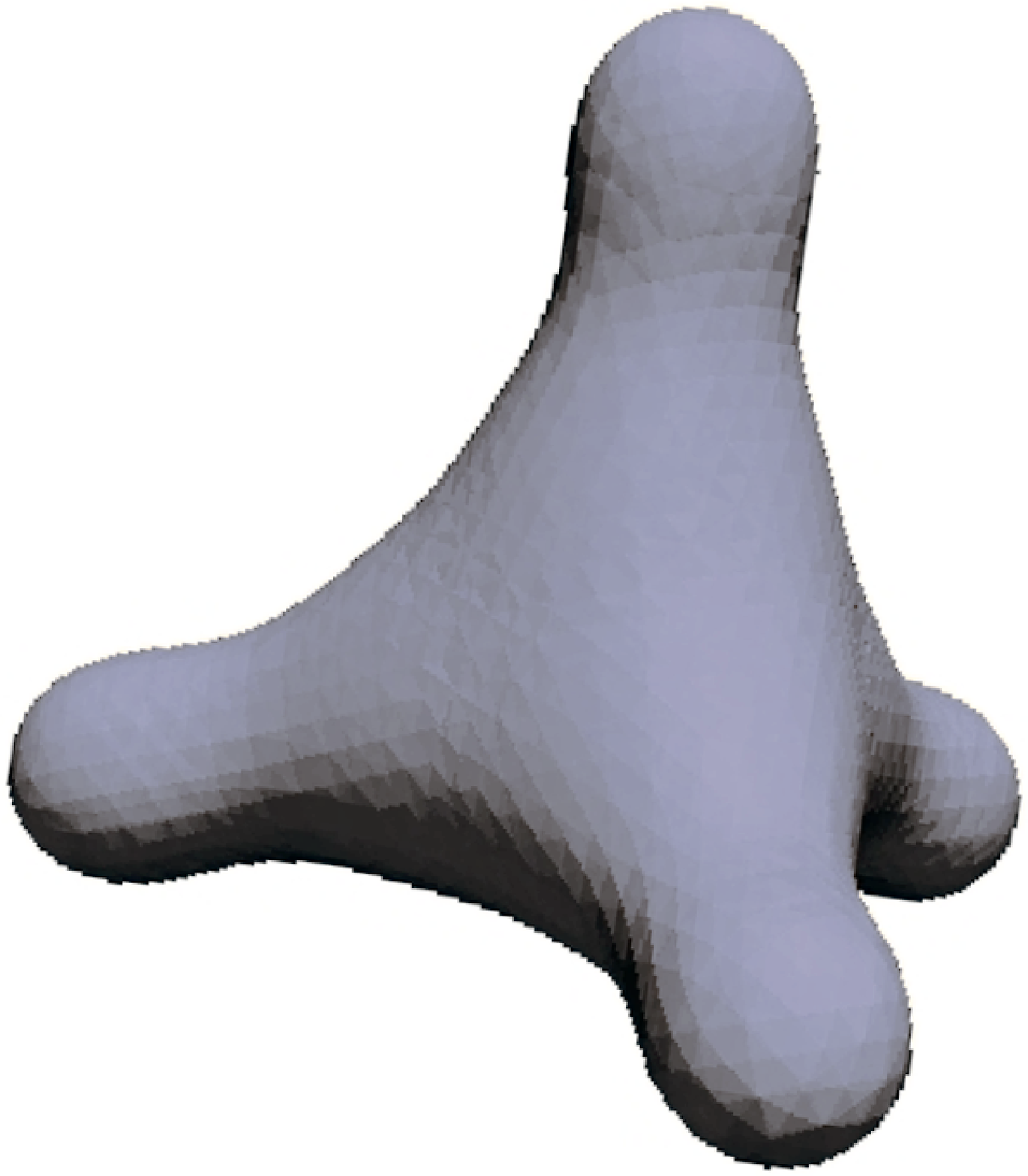}
 \includegraphics[width=5.5cm]{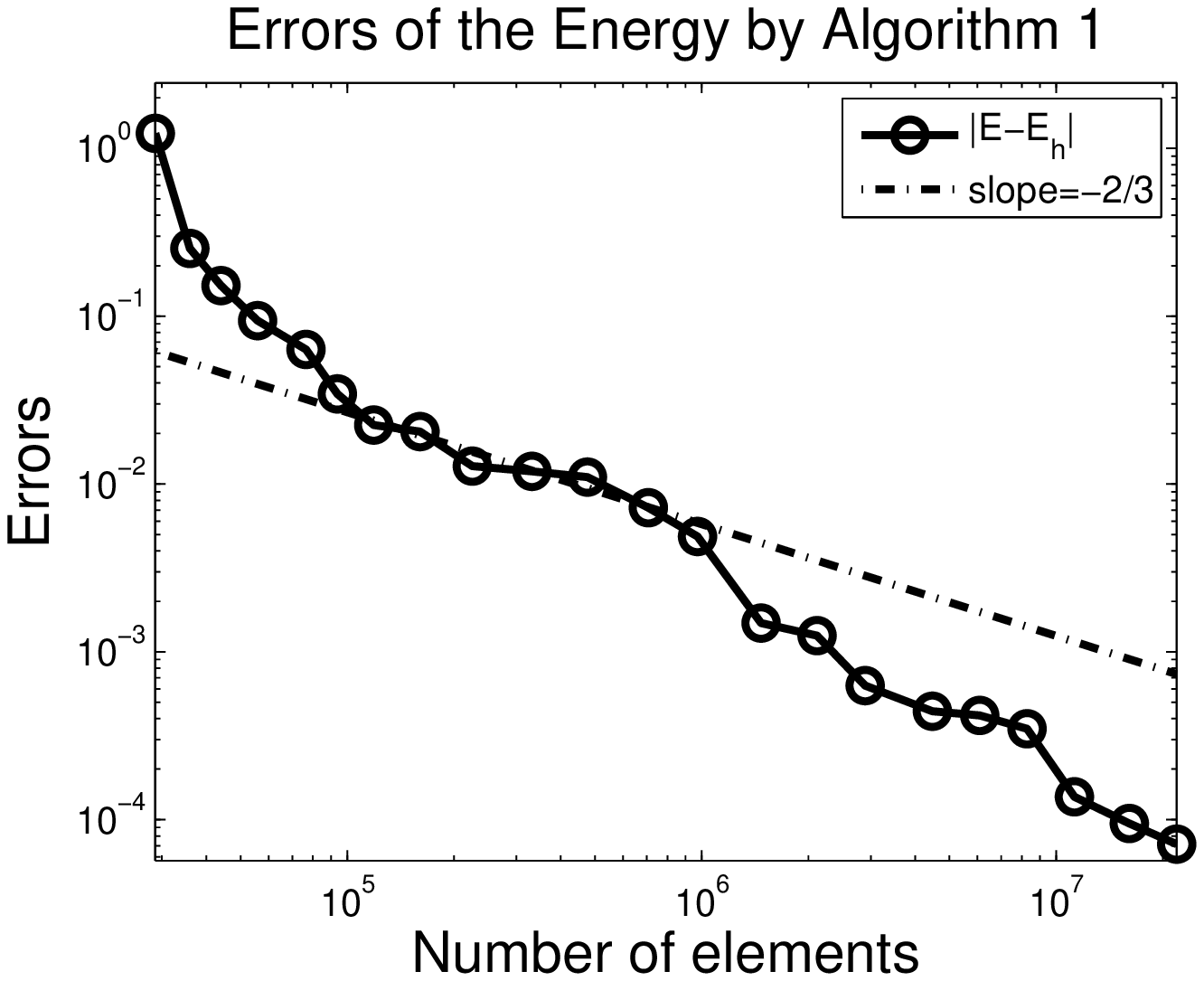}
\caption{Molecular structure of methane (left), the contour plot of the electron density (middle) and
errors of energy approximations (right)  of Algorithm \ref{Parallel_Aug_Subspace_Method} for Methane.}
\label{ex2-f1}
\end{figure}

\begin{figure}[htbp]
\centering
\includegraphics[width=4.6cm]{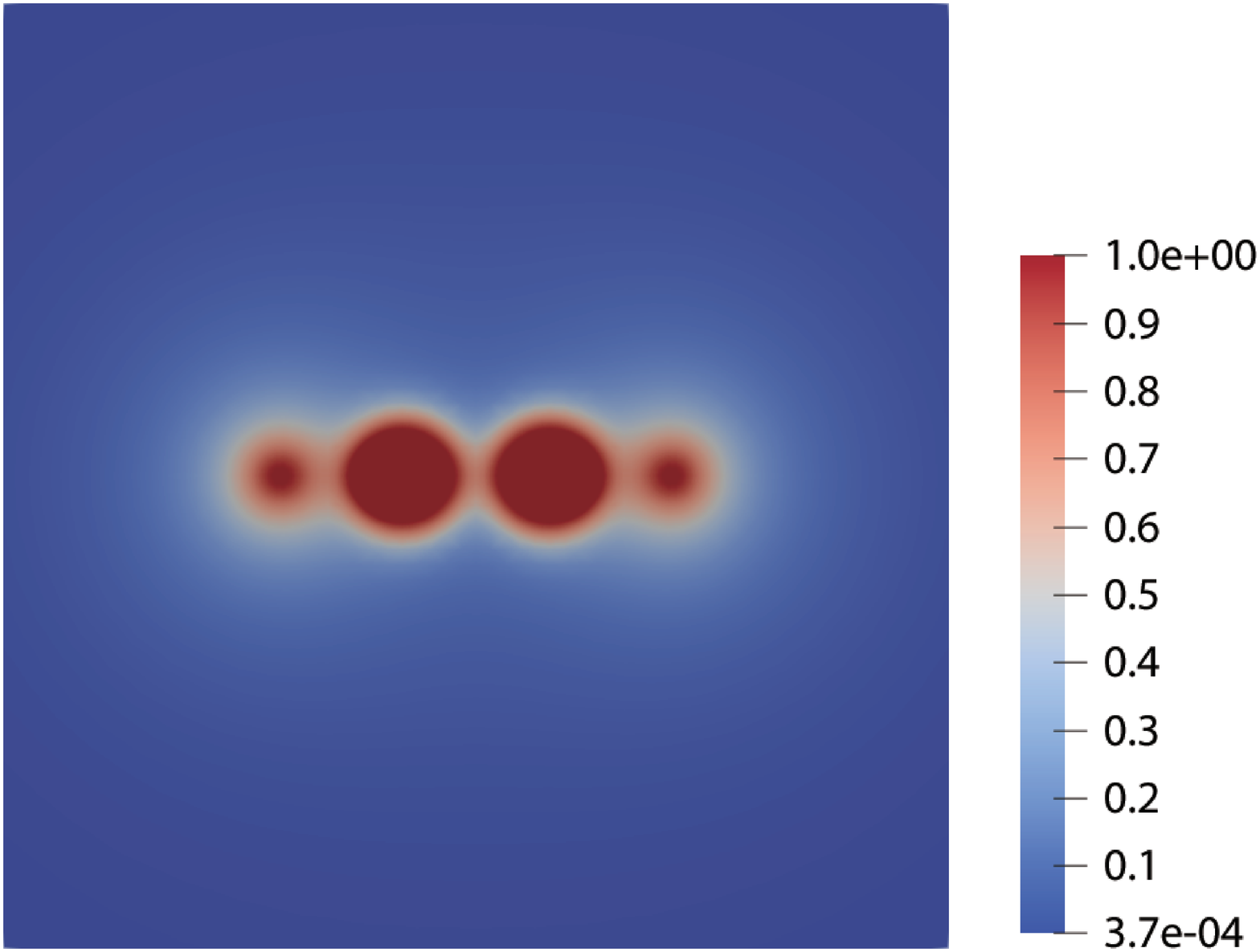}
\includegraphics[width=3.6cm]{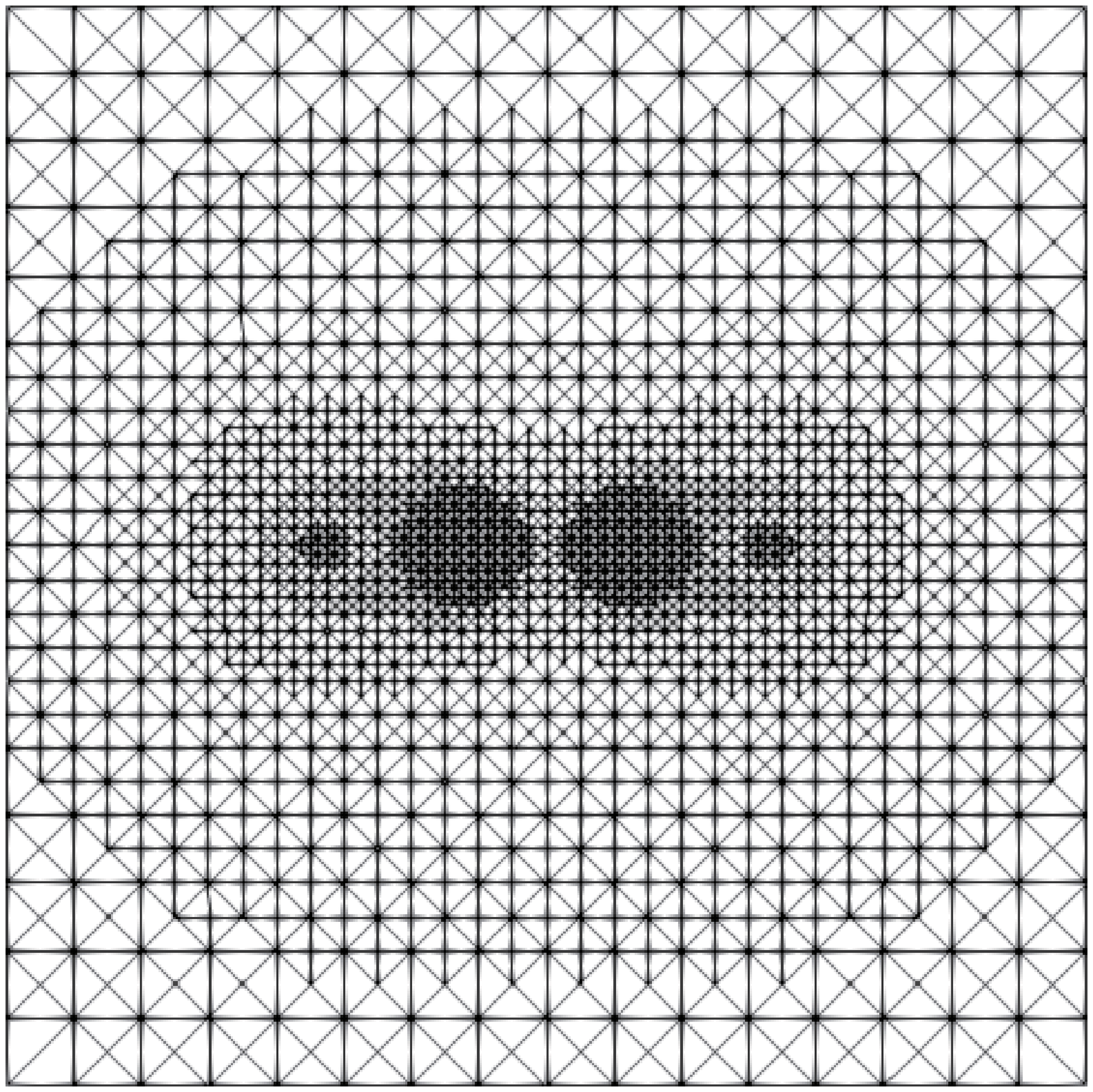}
\includegraphics[width=5.5cm]{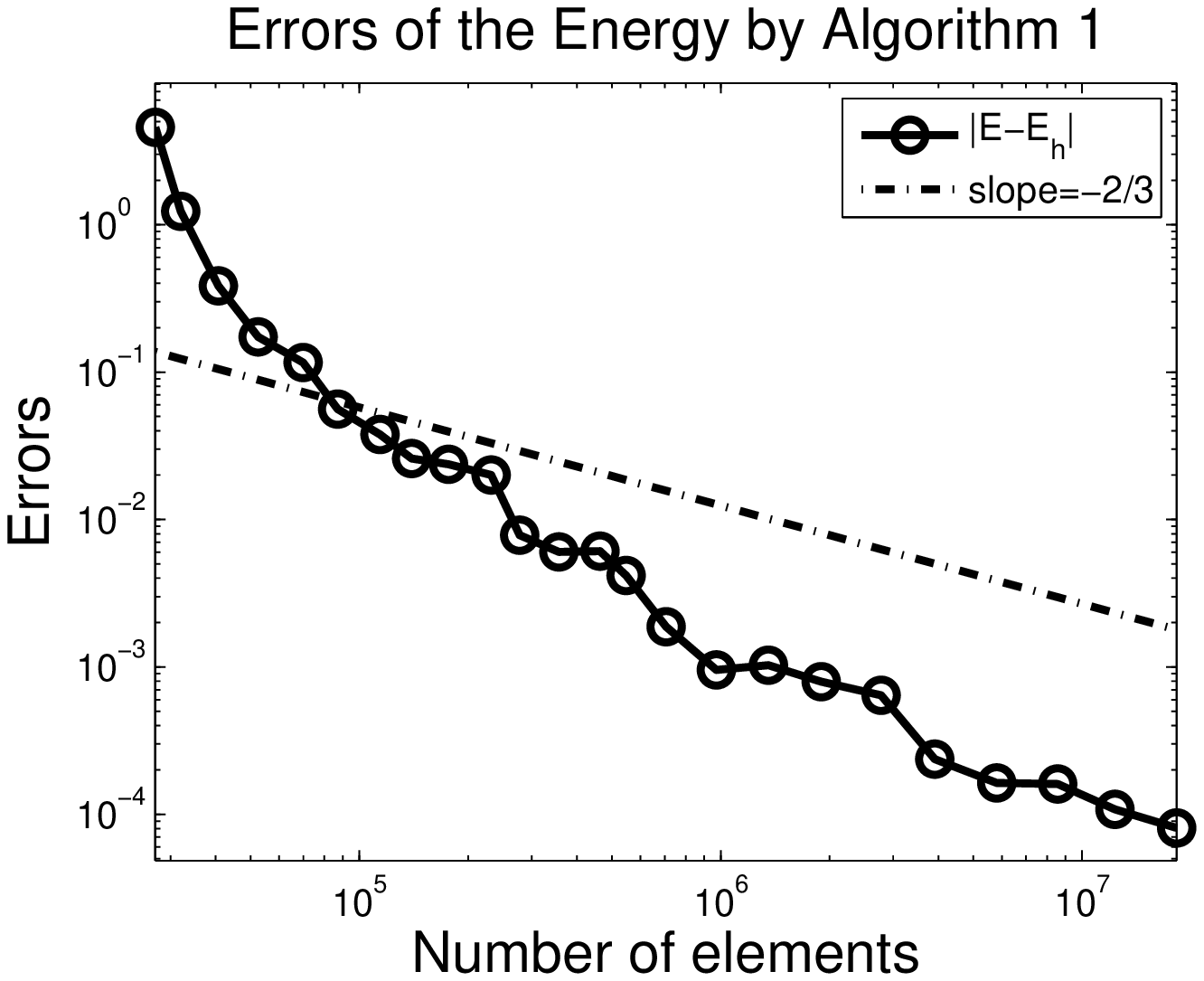}
\caption{ Contour plot of the electron density (left), the mesh after 15 adaptive refinements (middle) and
errors of energy approximations (right) of Algorithm \ref{Parallel_Aug_Subspace_Method} for Acetylene.}
\label{ex3-f1}
\end{figure}

\begin{figure}[htbp]
\centering
\includegraphics[width=4.8cm]{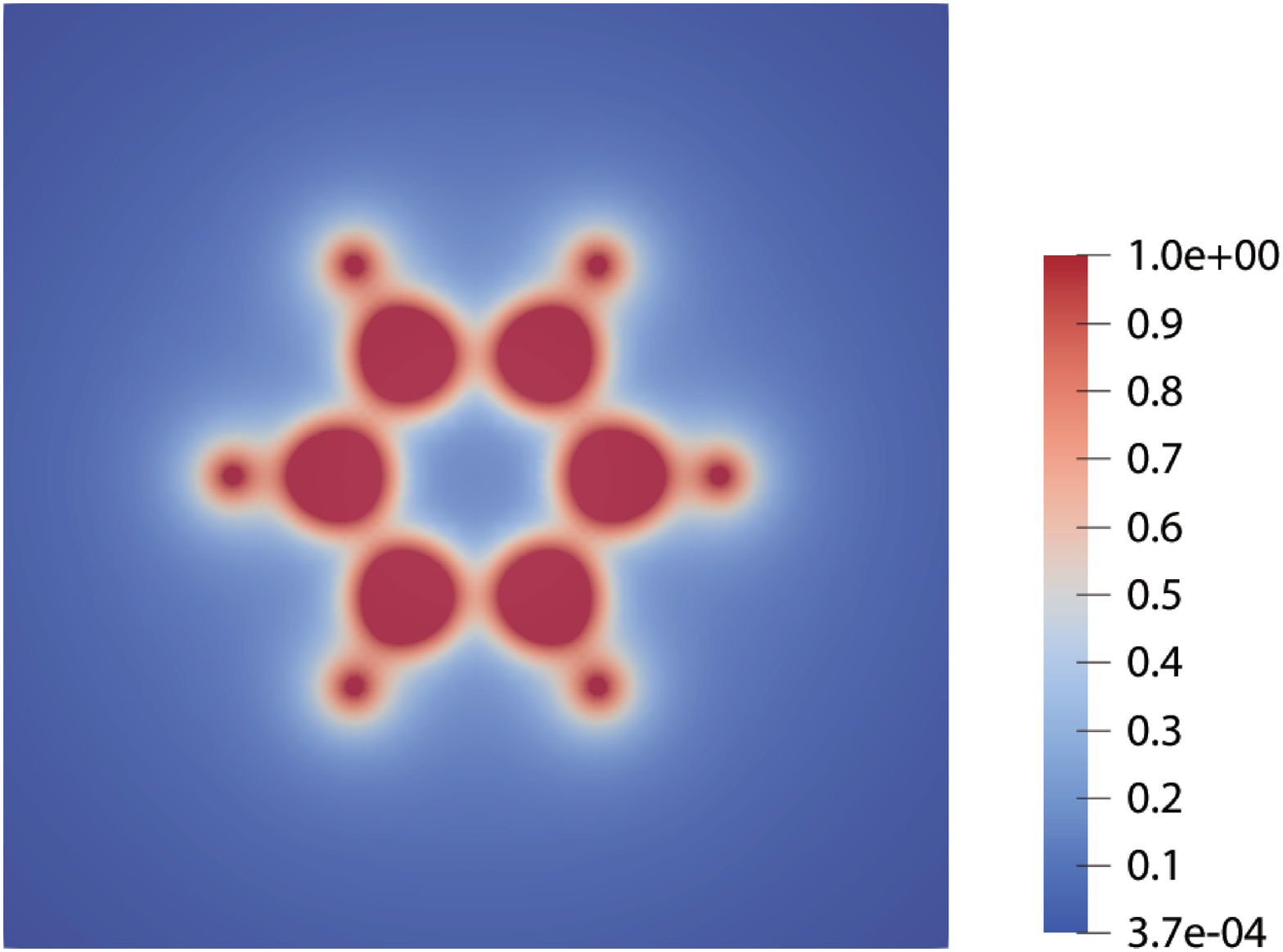}
\includegraphics[width=3.8cm]{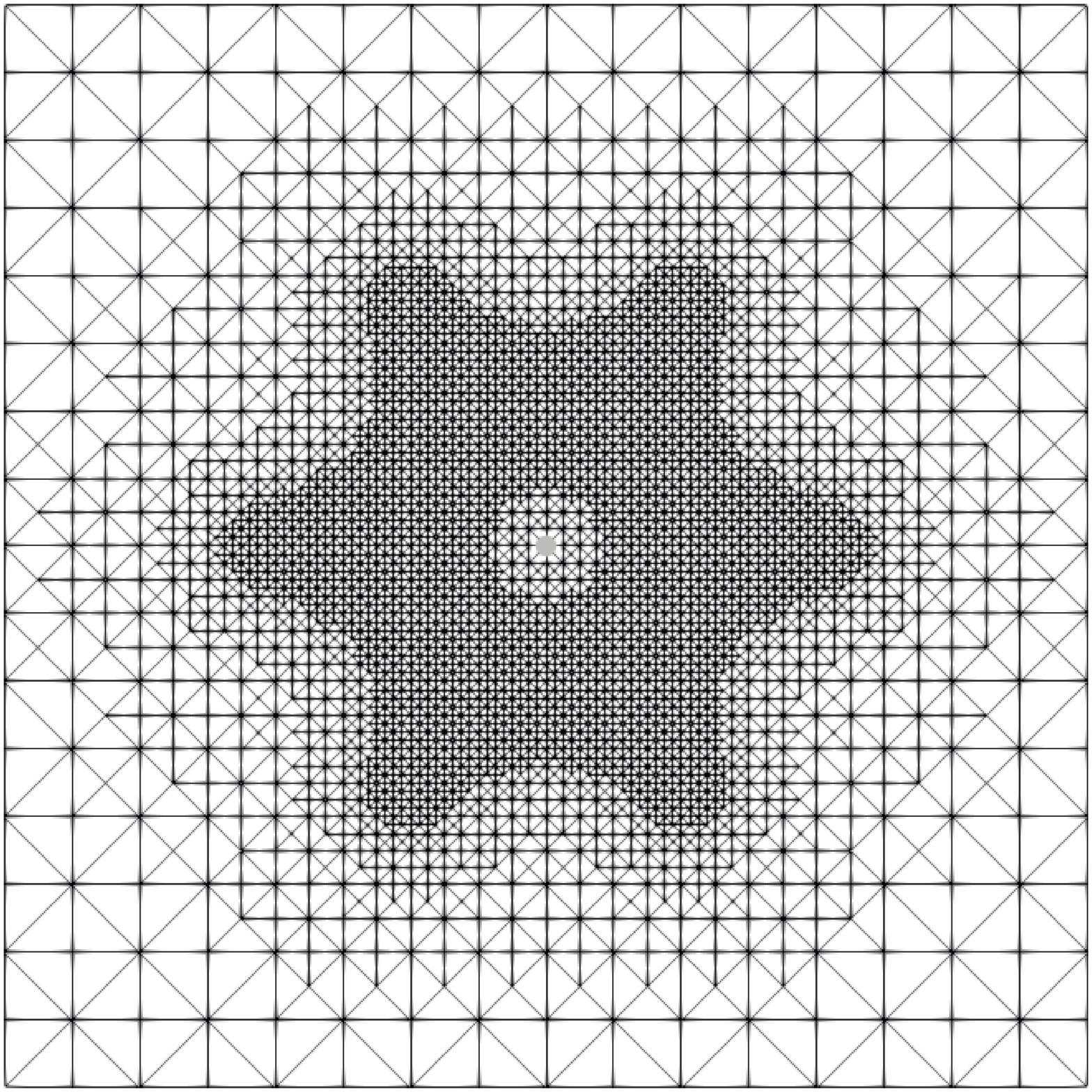}
\includegraphics[width=5.5cm]{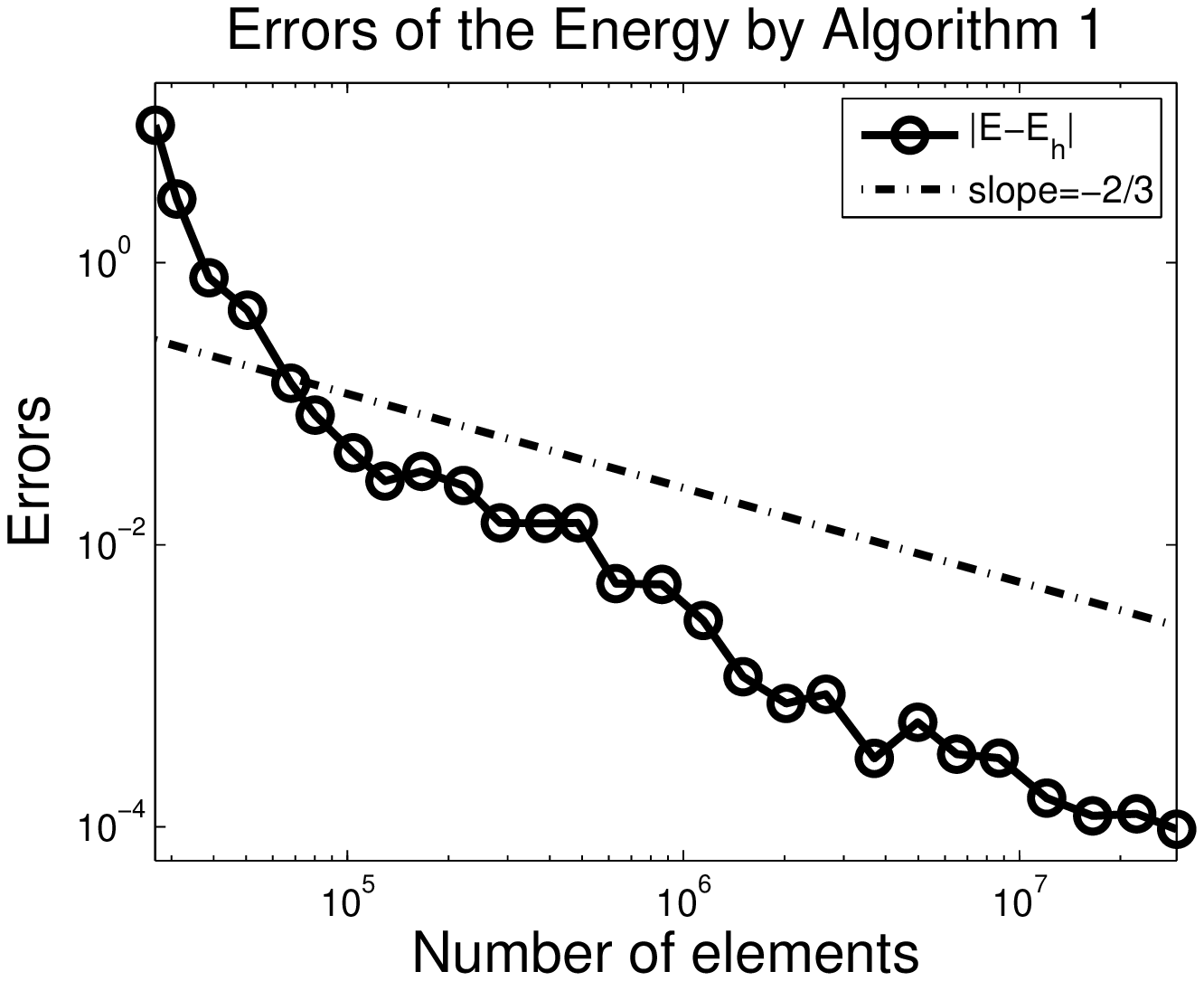}
\caption{ Contour plot of the electron density (left), the mesh after $15$ adaptive refinements (middle)
and errors of energy approximations  of Algorithm  \ref{Parallel_Aug_Subspace_Method} for Benzene.}
\label{ex4-f1}
\end{figure}

\begin{figure}[htbp]
\centering
\includegraphics[width=3.7cm]{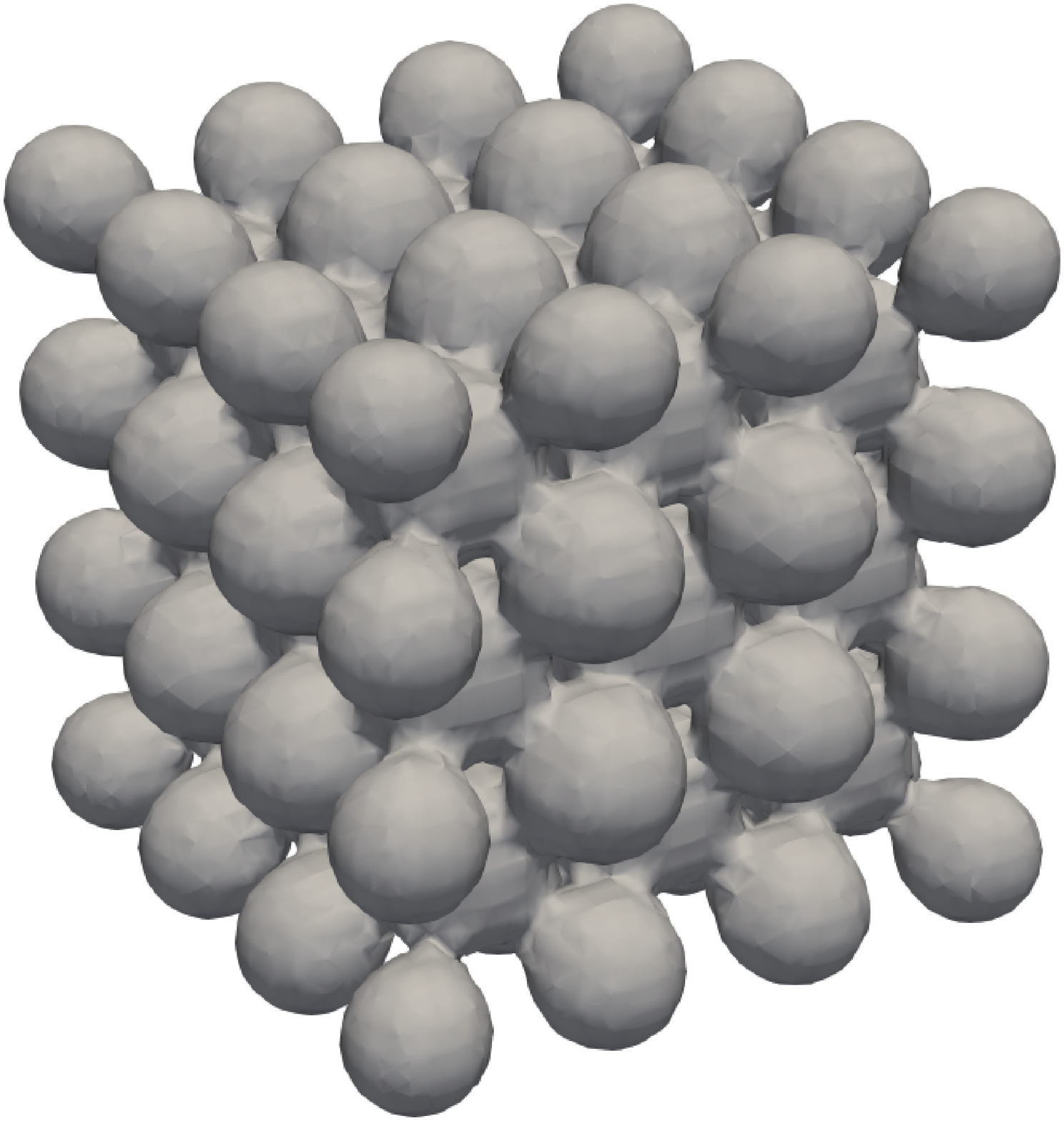} \ \ \ \ \
\includegraphics[width=5.5cm]{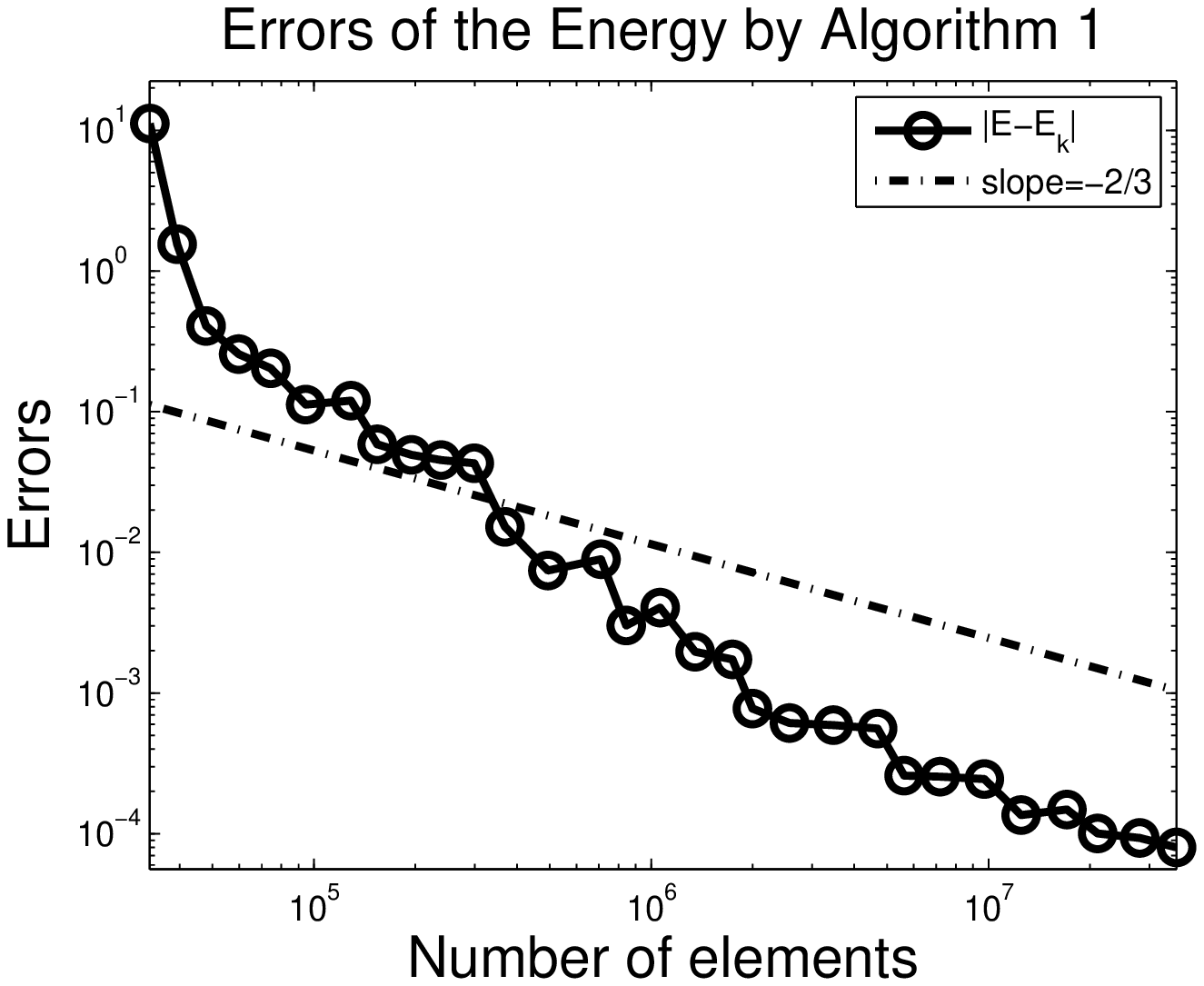}
\caption{Contour plot of the electron density (left) and errors of energy approximations (right) of Algorithm \ref{Parallel_Aug_Subspace_Method} for Sodium crystal.}
\label{ex5-f1}
\end{figure}

Figures \ref{ex1-f1},  \ref{ex2-f1}, \ref{ex3-f1}, \ref{ex4-f1} and \ref{ex5-f1}
gives the graphic exhibition of the numerical results by Algorithm
\ref{Parallel_Aug_Subspace_Method}.
The corresponding results associated with error estimates are presented in Figures
\ref{ex1-f2},  \ref{ex2-f1},  \ref{ex3-f1},  \ref{ex4-f1}  and \ref{ex5-f1}, which show that Algorithm
\ref{Parallel_Aug_Subspace_Method} is capable of deriving the optimal
error estimates.

\begin{table}[htbp]
\begin{center}
\begin{tabular}{|c|c|c|c|c|c|}\hline
Mesh Level &1& 2 &3 &4&5  \\ \hline
Inner Product &0.0000e-0& 1.9498e-5 & 1.0040e-5 &6.1033e-6&1.4316e-5 \\ \hline
\end{tabular}
\begin{tabular}{|c|c|c|c|c|c|}\hline
Mesh Level &6& 7 &8 &9&10  \\ \hline
Inner Product &3.2717e-6& 2.8455e-6 & 8.3587e-7 &1.9357e-6& 6.5374e-7 \\ \hline
\end{tabular}
\begin{tabular}{|c|c|c|c|c|c|}\hline
Mesh Level &11& 12 &13 &14&15  \\ \hline
Inner Product &1.7357e-6& 4.7698e-7 & 3.1525e-7 &2.2424e-7&2.8045e-7 \\ \hline
\end{tabular}
\begin{tabular}{|c|c|c|c|c|c|}\hline
Mesh Level &16&17&18 &19 &20  \\ \hline
Inner Product &8.8124e-7& 1.2938e-7 & 5.5421e-8 &2.5063e-7&4.8012e-8 \\ \hline
\end{tabular}
\begin{tabular}{|c|c|c|c|c|c|}\hline
Mesh Level &21&22&23 &24 &25  \\ \hline
Inner Product &3.7354e-8& 2.4355e-8 & 8.3322e-9&2.1567e-8 &5.3984e-9 \\ \hline
\end{tabular}
\end{center}
\caption{The inner products of the eigenfunctions corresponding to different eigenvalues for Hydrogen-Lithium.}
\label{ex1-table4}
\end{table}

\begin{table}[htbp]
\begin{center}
\begin{tabular}{|c|c|c|c|c|c|}\hline
Mesh Level &1& 2 &3 &4&5  \\ \hline
Inner Product &0.0000e-0& 3.5339e-5 & 2.6196e-5 &7.8513e-6&2.1434e-5 \\ \hline
\end{tabular}
\begin{tabular}{|c|c|c|c|c|c|}\hline
Mesh Level &6& 7 &8 &9&10  \\ \hline
Inner Product &5.4062e-6& 1.6346e-5 & 4.6419e-6 &2.7387e-6&2.6765e-6 \\ \hline
\end{tabular}
\begin{tabular}{|c|c|c|c|c|c|}\hline
Mesh Level &11& 12 &13 &14&15  \\ \hline
Inner Product &1.2886e-6& 1.1922e-6 & 8.6543e-7 &1.5147e-6&7.7643e-7 \\ \hline
\end{tabular}
\begin{tabular}{|c|c|c|c|c|c|}\hline
Mesh Level &16&17&18 &19 &20  \\ \hline
Inner Product &6.2544e-7& 3.3645e-7 & 1.9469e-7 &2.4036e-7&1.0117e-7 \\ \hline
\end{tabular}
\begin{tabular}{|c|c|c|c|c|c|}\hline
Mesh Level &21&22&23 &24 &25  \\ \hline
Inner Product &8.8124e-8& 1.8111e-7 & 2.5065e-8 &5.0711e-8&1.1934e-8 \\ \hline
\end{tabular}
\end{center}
\caption{The inner products of the eigenfunctions corresponding to different eigenvalues for Methane.}
\label{ex2-table2}
\end{table}

\begin{table}[htbp]
\begin{center}
\begin{tabular}{|c|c|c|c|c|c|}\hline
Mesh Level &1& 2 &3 &4&5  \\ \hline
Inner Product & 0.0000e-0 & 7.4681e-5 & 6.1964e-5 & 4.0481e-5 & 2.1546e-5 \\ \hline
\end{tabular}
\begin{tabular}{|c|c|c|c|c|c|}\hline
Mesh Level &6& 7 &8 &9&10  \\ \hline
Inner Product &5.4252e-6& 1.4153e-5 & 3.5523e-6 &1.0320e-5&4.1546e-6 \\ \hline
\end{tabular}
\begin{tabular}{|c|c|c|c|c|c|}\hline
Mesh Level &11& 12 &13 &14&15  \\ \hline
Inner Product &3.0646e-6& 2.5461e-6 & 9.4651e-7 &1.0313e-6&1.4862e-6 \\ \hline
\end{tabular}
\begin{tabular}{|c|c|c|c|c|c|}\hline
Mesh Level &16&17&18 &19 &20  \\ \hline
Inner Product &8.6463e-7& 6.5663e-7 & 7.6741e-7 &7.2418e-7&5.1564e-7 \\ \hline
\end{tabular}
\begin{tabular}{|c|c|c|c|c|c|}\hline
Mesh Level &21&22&23 &24 &25  \\ \hline
Inner Product & 8.5362e-8 & 5.2283e-8 &2.1813e-7&4.7199e-8&2.5502e-8 \\ \hline
\end{tabular}
\end{center}
\caption{The inner products of the eigenfunctions corresponding to different eigenvalues for Acetylene.}
\label{ex3-table3}
\end{table}

\begin{table}[htbp]
\begin{center}
\begin{tabular}{|c|c|c|c|c|c|}\hline
Mesh Level &1& 2 &3 &4&5  \\ \hline
Inner Product &0.0000e-0& 3.3343e-4 & 1.6784e-5 &8.3429e-4 & 7.7229e-5 \\ \hline
\end{tabular}
\begin{tabular}{|c|c|c|c|c|c|}\hline
Mesh Level &6& 7 &8 &9&10  \\ \hline
Inner Product &5.5827e-5& 2.4008e-5 & 2.7386e-5 &1.4684e-5& 8.1615e-6 \\ \hline
\end{tabular}
\begin{tabular}{|c|c|c|c|c|c|}\hline
Mesh Level &11& 12 &13 &14&15  \\ \hline
Inner Product &1.6854e-5& 5.3768e-6 & 3.6524e-6 &3.7026e-6&2.5145e-6 \\ \hline
\end{tabular}
\begin{tabular}{|c|c|c|c|c|c|}\hline
Mesh Level &16&17&18 &19 &20  \\ \hline
Inner Product & 2.4746e-6& 1.9326e-6 & 8.4358e-7 & 1.1354e-6& 7.1456e-7 \\ \hline
\end{tabular}
\begin{tabular}{|c|c|c|c|c|c|}\hline
Mesh Level &21&22&23 &24 &25  \\ \hline
Inner Product &1.1254e-6&4.6510e-7  & 5.9616e-7 &2.1382e-7&1.4798e-7 \\ \hline
\end{tabular}
\begin{tabular}{|c|c|c|c|c|c|}\hline
Mesh Level &26&27&28 &29 &30  \\ \hline
Inner Product &7.1642e-8& 1.3617e-7 & 8.4178e-8 &7.5464e-8&6.3578e-8 \\ \hline
\end{tabular}
\end{center}
\caption{The inner products of the eigenfunctions corresponding to different eigenvalues for Benzene.}
\label{ex4-table2}
\end{table}

\begin{table}[htbp]
\begin{center}
\begin{tabular}{|c|c|c|c|c|c|}\hline
Mesh Level &1& 2 &3 &4&5  \\ \hline
Inner Product &0.0000e-0& 8.4713e-4 & 8.0158e-4 &6.0454e-4&4.1687e-4 \\ \hline
\end{tabular}
\begin{tabular}{|c|c|c|c|c|c|}\hline
Mesh Level &6& 7 &8 &9&10  \\ \hline
Inner Product &8.5423e-5& 2.3185e-4 & 6.3291e-5 &4.8652e-5&3.5614e-5 \\ \hline
\end{tabular}
\begin{tabular}{|c|c|c|c|c|c|}\hline
Mesh Level &11& 12 &13 &14&15  \\ \hline
Inner Product &4.3523e-5& 2.5875e-5 & 1.7135e-5 &8.9253e-6&7.9633e-6 \\ \hline
\end{tabular}
\begin{tabular}{|c|c|c|c|c|c|}\hline
Mesh Level &16&17&18 &19 &20  \\ \hline
Inner Product & 1.5853e-5& 1.1443e-5 & 5.6537e-6 &3.2123e-6&3.7452e-6 \\ \hline
\end{tabular}
\begin{tabular}{|c|c|c|c|c|c|}\hline
Mesh Level &21&22&23 &24 &25  \\ \hline
Inner Product & 2.2873e-6 & 1.1564e-6 &  7.5720e-7 & 1.4354e-6 & 6.5264e-7 \\ \hline
\end{tabular}
\begin{tabular}{|c|c|c|c|c|c|}\hline
Mesh Level &26&27&28 &29 &30  \\ \hline
Inner Product &5.1215e-7& 3.8238e-7 & 3.9654e-7 &2.1426e-7&1.0964e-7 \\ \hline
\end{tabular}
\end{center}
\caption{The inner products of the eigenfunctions corresponding to different eigenvalues  for Sodium crystal.}
\label{ex5-table2}
\end{table}

Finally, in order to show the efficiency of Algorithm
\ref{Parallel_Aug_Subspace_Method} more clearly, we compare its
computational time with that of the standard AFEM and the standard
multilevel correction adaptive method developed in \cite{HuXieXu}.  In
our test, we set the same accuracy of energy for all the adopted
algorithms.
The CPU time (in seconds) is provided in Tables
\ref{timetable1} and \ref{timetable2}.  From Tables \ref{timetable1}
and \ref{timetable2}, we can find that Algorithm
\ref{Parallel_Aug_Subspace_Method} and the standard multilevel
correction adaptive algorithm all have better efficiency over the
standard AFEM.  The saving of the computing time comes from that
solving large-scale nonlinear eigenvalue problems is avoided in this
two algorithms.  Besides, we also can find that Algorithm
\ref{Parallel_Aug_Subspace_Method} has a large advantage over the
standard multilevel correction adaptive algorithm based on the new
computing strategy developed in this paper.

Furthermore, comparing the CPU time of different models presented in
Tables \ref{timetable1} and \ref{timetable2}, we can find the
advantage of Algorithm \ref{Parallel_Aug_Subspace_Method} is more
obvious for more complicated models. This is because a more
complicated model needs more SCF iteration numbers and this will takes
a significant amount of time for the classical methods; while for
Algorithm \ref{Parallel_Aug_Subspace_Method}, this can be solved
efficiently by using the separately handling method for the nonlinear terms and the efficient implementing schemes for the inner iterations.
In addition, comparing Tables \ref{timetable1} and \ref{timetable2}, we can find the
advantage of Algorithm \ref{Parallel_Aug_Subspace_Method} becomes more
obvious when the accuracy of energy improves.
This is because the computing time for solving the small-scale nonlinear eigenvalue problem (\ref{Nonlinear_Eig_Hh_Full}) is
gradually negligible along with the refinement of mesh.

\begin{table}[htbp]
\begin{center}
\begin{tabular}{|c|c|c|c|}\hline
\diagbox{Atom}{Time} & Time of standard AFEM & \tabincell{c}{Time of standard multilevel  \\ correction adaptive algorithm}
& Time of Algorithm 1 \\ \hline
Hydrogen-Lithium  &188.93   &  95.74   & 94.67   \\ \hline
Methane       &  738.85 &    321.13    &  272.31\\ \hline
Acetylene       & 3426.61  &   685.23     &  284.78\\ \hline
Benzene       &  7311.49 &   1354.75     & 646.16 \\ \hline
Sodium crystal & 17351.72  &   2629.95     &  1069.04\\ \hline
\end{tabular}
\end{center}
\caption{The computational time (in seconds) of the standard AFEM, the standard multilevel correction adaptive algorithm and
Algorithm \ref{Parallel_Aug_Subspace_Method} with the same energy accuracy 5E-3.}
\label{timetable1}
\end{table}

\begin{table}[htbp]
\begin{center}
\begin{tabular}{|c|c|c|c|}\hline
\diagbox{Atom}{Time} & Time of standard AFEM & \tabincell{c}{Time of standard multilevel \\ correction adaptive algorithm}
& Time of Algorithm 1 \\ \hline
Hydrogen-Lithium  &2986.38   &  605.29   & 508.39   \\ \hline
Methane       &  16052.52 &     2494.06    &  1648.14\\ \hline
Acetylene       & 28091.26 &    3942.35     &  1736.56\\ \hline
Benzene       &  142339.41 &      12374.82     &  3821.42\\ \hline
Sodium crystal &  - &     26037.58     &  6562.57\\ \hline
\end{tabular}
\end{center}
\caption{The computational time (in seconds) of the standard AFEM, the standard multilevel correction adaptive algorithm and
Algorithm \ref{Parallel_Aug_Subspace_Method} with the same energy accuracy 1E-4. The symbol ``-" means the computer runs out of memory.}
\label{timetable2}
\end{table}

\section{Conclusion}
In this paper, we propose an accelerating multilevel correction adaptive
finite element method for the Kohn-Sham equation.
The new algorithm benefits from two acceleration strategies.
The first one is to separately handle the nonlinear Hartree potential and exchange-correlation potential,
which can be solved efficiently by outer iteration and inner iteration. respectively.
The second one is to parallelize the algorithm in an eigenpairwise approach.
Compared with previous results, a significant improvement of numerical efficiency can be observed from plenty of numerical
experiments, which make the new method more suitable for the practical problems.



\end{document}